\documentclass{scrartcl}
\usepackage{amsfonts}
\usepackage{amsmath}
\usepackage{amssymb}
\usepackage{hyperref}
\usepackage{cleveref}
\usepackage{mathrsfs}
\usepackage{bbm}
\usepackage{mathtools}
\usepackage{enumitem}
\usepackage{tikz}
\usetikzlibrary{arrows}
\usepackage{pgfplots}
\usepackage{wrapfig}
\usepackage{float}
\allowdisplaybreaks

\newcommand{\bbbr}{\mathbb{R}}
\newcommand{\bbbn}{\mathbb{N}}
\newcommand{\bbbz}{\mathbb{Z}}

\newcommand{\Ngen}{\mathcal{N}}

\newcommand{\htwo}{\mathcal{H}^{2}}
\newcommand{\ctI}{T_{I}}

\newcommand{\ctJ}{T_{J}}

\newcommand{\ctIJ}{T_{I \times J}}
\newcommand{\that}{\widehat{T}}
\newcommand{\ahat}{\widehat{A}}
\newcommand{\const}{\mathsf{C}}
\newcommand{\smallconst}{\mathsf{c}}
\newcommand{\interpol}{\mathfrak{I}}
\newcommand{\lagrange}{\mathfrak{L}}
\newcommand{\Root}{\mathscr{R}}
\newcommand{\Leaves}{\mathscr{L}}
\newcommand{\Bsf}{\mathsf{B}}
\newcommand{\diam}{\mathop{\operatorname{diam}}\nolimits}
\newcommand{\dist}{\mathop{\operatorname{dist}}\nolimits}
\newcommand{\supp}{\mathop{\operatorname{supp}}\nolimits}

\newcommand{\level}{\mathop{\operatorname{level}}\nolimits}
\newcommand{\sons}{\mathop{\operatorname{sons}}\nolimits}

\newcommand{\true}{\mathop{\operatorname{true}}\nolimits}
\newcommand{\false}{\mathop{\operatorname{false}}\nolimits}
\newtheorem{theorem}{Theorem}
\newtheorem{lemma}{Lemma}
\newtheorem{definition}{Definition}
\newtheorem{remark}{Remark}
\newtheorem{assumption}{Assumption}
\title{\texorpdfstring{$\mathcal{H}^2$}{H2}-matrices for
       translation-invariant kernel functions}

\author{Steffen B\"orm and Janne Henningsen}

\begin{document}
	
\maketitle

\begin{abstract}
Boundary element methods for elliptic partial differential equations
typically lead to boundary integral operators with translation-invariant
kernel functions.
Taking advantage of this property is not straightforward if general
unstructured meshes and general basis functions are used, since we need the
supports of these basis functions to be contained in a hierarchy of
subdomains with translational symmetry.

In this article, we present a modified construction for
$\mathcal{H}^2$-matrices on unstructured quasi-uniform meshes that uses
translation-invariance to significantly reduce the storage requirements
for the farfield representation.

We construct a nested hierarchy of axis-parallel boxes so that
translational symmetry is preserved and prove optimal-order complexity
estimates under moderate assumptions. In particular, we need only one weak assumption for proving that the entire farfield requires only $ \mathcal{O}(\log(n)) $ coefficients.

It should be mentioned that, since we are working with an unstructured
mesh and general basis functions, the nearfield of the matrix
still requires $\mathcal{O}(n)$ units of storage.
\end{abstract}

This work was funded by the DFG in project BO 3289/7-1.



\section{Introduction}
\label{se:introduction}

Boundary element methods are an attractive technique for handling
homogeneous linear partial differential equations, e.g., the
Laplace, Lam\'e, Helmholtz, or Maxwell equations
\cite{HA95,HSWE08,SASC11}, on a domain
$\Omega \subseteq \bbbr^d$ with $ d \in \{2,3\}$.
A Galerkin discretization of related boundary integral operators
based on test functions $(\varphi_i)_{i\in I}$ and trial functions
$(\psi_j)_{j\in J}$ leads to matrices $G\in\bbbr^{I \times J}$ given by
\begin{align}\label{eq:matrix_entries}
  G_{ij} &:= \int_\Gamma \varphi_i(x)
            \int_\Gamma g(x,y) \psi_j(y) \,dy\,dx &
  &\text{ for all } i \in I,\ j \in J,
\end{align}
where $\Gamma := \partial\Omega$ denotes the boundary of 
$\Omega$ and $g\colon\bbbr^{d}\times\bbbr^{d}\to\bbbr$ is a
kernel function.

In typical applications, the kernel function $g$ is non-zero
almost everywhere, and the matrix $G$ is therefore densely populated.
In order to reduce the storage requirements and computational
complexity, compression techniques are employed, e.g., the
panel-clustering technique \cite{HANO89}, the fast multipole
method \cite{RO85,GRRO87,AN92,GRRO97,BIYIZO04},
interpolation \cite{GI01,BOHA02a,BOLOME02,BOME15,DAFO09},
algebraic approximations \cite{TY96,TY99,BE00a,BEVE12,CORAZO15},
or hybrid methods \cite{GIRO02,BOGR04,BOCH14}.

All of these techniques split the matrix $G$ into submatrices
that can belong either to the \emph{nearfield}
or the \emph{farfield}.
Nearfield submatrices are small and can be stored directly, while
farfield submatrices can be large and have to be approximated, e.g.,
by low-rank matrices that can be stored efficiently in factorized
form.
Unfortunately, even the factorized form still requires a large
amount of storage, and while recompression techniques
\cite{BOHA03,GR04,BOBO18} can help, they require additional
computational work.

In this article, we pursue an alternative approach frequently used in the context of
fast multipole methods for particle systems:
if the kernel function $g$ is invariant under translation,
i.e., if
\begin{align}\label{eq:translation}
  g(x,y) &= g(x+c,y+c) &
  &\text{ for all } x,y,c\in\bbbr^d,
\end{align}
we can modify interpolation and certain hybrid methods in a
way that reduces the storage requirements for the farfield matrices
to $\mathcal{O}(k^2 \log(n))$, where $k$ denotes the rank of the
approximation and $n:=\max\{|I|,|J|\}$ the maximal dimension
of the matrix $G$.

This task is straightforward for particle systems with approximately
uniformly distributed particles, since each particle is represented
by just one point in space, but it is significantly more challenging
for Galerkin discretizations, since every matrix entry depends on
the supports of the basis functions $\varphi_i$ and $\psi_j$, and the
supports of different basis functions may overlap or have a non-trivial shape.
In some implementations, this problem is circumvented by replacing
the integrals in \cref{eq:matrix_entries} by quadrature, thus reducing
the computation to a sum of kernel evaluations closely related to
particle methods \cite{BESC21}.
While this is certainly an elegant approach and allows BEM codes
to take advantage of highly sophisticated implementations of the fast
multipole method, the number of quadrature points may grow excessively
large if higher quadrature orders are used to keep up with the
discretization error as the grid is refined.

In this article, we present a different approach:
Initially, we construct a hierarchy of axis-parallel boxes that cover representative points
of the
supports of the basis functions, and we ensure that all
boxes on a given level are translations of the same reference box
with a fixed displacement step size. By suitably enlarging the reference boxes and applying the aforementioned translations to these enlarged reference boxes we then construct a corresponding hierarchy of axis-parallel boxes that cover the entire supports. 

If we apply standard techniques like interpolation \cite{GI01,BOHA02a} or
Green quadrature \cite{BOCH14} to this structure, matrices corresponding
to translation-equivalent pairs of boxes are identical and therefore
have to be stored only once.
As long as the enlargement of the reference boxes is sufficiently bounded, this allows us to reduce the corresponding storage requirements to $\mathcal{O}(k^2)$
for each level. No further assumption is needed for proving this result.

This approach also reduces the assembly time, since the matrices
have to be set up only once.
The run-time for the matrix-vector multiplication may benefit, too, since
fewer data has to be moved between main memory and the processor,
although the number of arithmetic operations usually grows since guaranteeing translational symmetry limits
our flexibility when constructing the boxes.

That is why a major part of this article is devoted to proving that
our proposed construction does not change the asymptotically optimal
complexity of the $\mathcal{H}^2$-matrix method.

In \cref{se:h2matrices} we introduce $\mathcal{H}^2$-matrices, a suitable
representation of our matrix approximation, and outline the modifications
required to take full advantage of translation-invariance.
\Cref{se:complexity} contains the key result of this article: we prove
that the entire farfield of the modified $\mathcal{H}^2$-matrices
requires only $\mathcal{O}(k^2 \log(n))$ coefficients if
the enlargement of the reference boxes is sufficiently bounded.
The nearfield matrices and the leaf matrices are directly connected
to the unstructured surface mesh and therefore cannot take advantage of
the translation-invariance property, so we store these matrices
explicitly using $\mathcal{O}(n k)$ units of storage.
\Cref{se:experiments} illustrates the advantage of the modified
representation in a series of numerical experiments.

\section{\texorpdfstring{$\mathcal{H}^2$-matrices for translation-invariant
    kernel functions}
  {H2-matrices for translation-invariant kernel functions}}
\label{se:h2matrices}

In order for our approximation technique to achieve the desired accuracy, the kernel function $ g $ needs to be \emph{asymptotically smooth}, i.e., there have to be constants $ \const_{as} \in \bbbr_{\geq0} $, $ \smallconst_{0} \in \bbbr_{>0} $ and a \emph{singularity degree} $ \sigma \in \bbbn $ such that for all $ (x,y) \in \bbbr^{d} \times \bbbr^{d} $  with $ x \neq y $ the following condition is fulfilled:
\begin{equation}\label{ineq:as}
\left|\partial^{\nu}_{\iota}g(x,y)\right| \leq 
\const_{as}\frac{(\sigma - 1 + \nu)!\,\smallconst_{0}^{\nu}}{\|x - y\|_{2}^{\nu+\sigma}}
\qquad      
\text{for all $ \nu \in \bbbn_{0}, \iota \in \{1,...,2d\}$}.
\end{equation}
Kernel functions occurring in typical applications are known to be asymptotically smooth, see for example \cite[Appendix E]{HA15} and \cite{BAHA05}. Additionally, we assume that $ g $ satisfies the translation-invariance property (\ref{eq:translation}).
\begin{definition}[Tree notations]\label{tree notations}
	Let $ T $ be a tree. We use the notation $ t \in T $ for "$ t $ is a node in $ T $".  For each node $ t \in T $ we denote the set of its sons by $ \sons_{T}(t) $. The set of leaves of $ T $ is denoted by
	$ \Leaves_{T} := \{t \in T : \sons_{T}(t) = \emptyset\}. $
	We denote the \emph{root} of $ T $ by $ \Root_{T} $ and define the
        \emph{level} of a node by
	\begin{equation*}
	\level_{T}(t) := 
	\begin{cases}
	0 &\text{if $ t = \Root_{T} $},\\
	\level_{T}(\check{t}) + 1 &\text{if $ t $ has a father $ \check{t} \in T$}
	\end{cases}
	\qquad \text{for all $ t \in T$}.
	\end{equation*}
	For all $ \ell \in \bbbn_{0} $ we define
	$  T^{(\ell)} := \{t \in T : \level_{T}(t) = \ell\}  $.
\end{definition}

The first phase of our approximation procedure relies on
\emph{characteristic points} $ (x_i)_{i\in I} $ and $ (y_{j})_{j \in J} $ satisfying
$ x_{i} \in \supp(\varphi_{i}) $ for all $ i \in I $ and
$ y_{j} \in \supp(\psi_{j}) $ for all $ j \in J $.
These points are split hierarchically into a hierarchy of boxes that
will then give rise to a decomposition of the matrix into submatrices.

We choose a maximal level $ \ell_{max} \in \bbbn_{\geq d}$ and construct trees
$\ctI$ and $\ctJ$ of closed axis-parallel boxes in $\bbbr^d$ satisfying
\begin{equation}\label{translational}
t \,=\, r^{(\ell)} + \underbrace{\delta^{(\ell)} \odot p_{t}}_{=:\, m_{t}}
\,=\, r^{(\ell)} + m_{t} \qquad \text{for all $ t \in \ctI^{(\ell)}\cup\ctJ^{(\ell)} $}
\end{equation}
on every level $ \ell \in \{0,...,\ell_{max}\} $, where $ (r^{(\ell)})_{\ell=0}^{\ell_{max}} $ is a family of \emph{reference boxes}, $ \delta^{(\ell)} \in \bbbr^{d}$ consists of the interval lengths of $ r^{(\ell)}$ for each $ \ell \in \{0,...,\ell_{max}\} $ and $ (p_{t})_{t\in \ctI \cup \ctJ} $ are suitable integer vectors.

We start with computing
an axis-parallel box $ r^{(0)} \subseteq \, \bbbr^{d}$ containing $ \Gamma $.
In order to create the isotropic boxes required by our complexity
analysis, we then compute $ (r^{(\ell)})_{\ell=1}^{\ell_{max}} $ recursively
by consecutive splitting in each coordinate direction.
The corresponding cyclic sequence of splitting directions is given by
\begin{align}\label{regular_splitting}
\iota^{(0)} := 1, \qquad \quad
\iota^{(\ell)}:=
\begin{cases}
\iota^{(\ell-1)} + 1 & \text{if $ \iota^{(\ell-1)} < d $},\\
1 & \text{if $ \iota^{(\ell-1)} = d $}
\end{cases}
\qquad \text{for all $ \ell \in \bbbn. $}
\end{align}
Assuming that
$ r^{(\ell)} = [a^{(\ell)}_{1}, b^{(\ell)}_{1}] \times ... \times [a^{(\ell)}_{d}, b^{(\ell)}_{d}] $
on a level $ \ell \in \{0,...,\ell_{max}-1\} $ is already given, we compute 
\begin{align}\label{halved_ref_box}
r^{(\ell+1)} 
:= \{z \in r^{(\ell)} \,:\, z_{\iota^{(\ell)}} \leq c^{(\ell)}_{\iota^{(\ell)}}\}
&& 
\text{with}
&&
c^{(\ell)}_{\iota^{(\ell)}} := \frac{a^{(\ell)}_{\iota^{(\ell)}} + b^{(\ell)}_{\iota^{(\ell)}}}{2}.
\end{align}
If $ \ell+1 < \ell_{max} $ holds, we proceed with $ r^{(\ell+1)} $ by
recursion. Otherwise, we stop.

Next, we create the tree $ \ctI $ by shifting the reference boxes
$r^{(\ell)}$ in a regular pattern.
In this stage, we also construct subsets $ (I_{t})_{t \in \ctI} $ of $ I $
such that
\begin{align}\label{index_subsets}
  x_i &\in t &\text{ for all } i\in I_t.
\end{align}
Beginning with $ r^{(0)} $ as the root $ \Root_{\ctI} $ of $ \ctI $, $ p_{r^{(0)}} := 0 $ and $ I_{r^{(0)}} := I$, we construct $ \ctI $ also recursively: Assuming that a box 
$ t = [a_{t,1}, b_{t,1}] \times ... \times [a_{t,d}, b_{t,d}] \in \ctI^{(\ell)} $ on a level $ \ell \in \{0,...,\ell_{max}-1\} $ with $ t = r^{(\ell)} + \delta^{(\ell)} \odot p_{t} $ for a vector $ p_{t} \in \bbbn_{0}^{d}$ 
and a corresponding non-empty subset $ I_{t} $ of $ I $ with
(\ref{index_subsets}) are already given, we compute the midpoint $ c_{t, \iota^{(\ell)}} := (a_{t,\iota^{(\ell)}} + b_{t,\iota^{(\ell)}})/2$ and
\begin{align*}
	t_{1} &\,:=\, \{z \in \bbbr^{d} : z_{\iota^{(\ell)}} \leq c_{t,\iota^{(\ell)}}\} \,=\, r^{(\ell+1)} + \delta^{(\ell+1)} \odot p_{t_{1}}, \\
	t_{2} &\,:=\, \{z \in \bbbr^{d} : z_{\iota^{(\ell)}} \geq c_{t,\iota^{(\ell)}}\} 
	\,=\, r^{(\ell+1)} + \delta^{(\ell+1)} \odot p_{t_{2}}
\end{align*}
with $ p_{t_{1}}, p_{t_{2}} \in \bbbn_{0}^{d} $ given by
\begin{equation}\label{displacement_sons}
p_{t_{1}, \iota} := 
\begin{cases}
2p_{t,\iota} & \text{if $ \iota = \iota^{(\ell)} $} \\
p_{t,\iota}  & \text{if $ \iota \neq \iota^{(\ell)} $}
\end{cases},
\quad
p_{t_{2}, \iota} := 
\begin{cases}
2p_{t,\iota} + 1 & \text{if $ \iota = \iota^{(\ell)} $} \\
p_{t,\iota}  & \text{if $ \iota \neq \iota^{(\ell)} $}
\end{cases}
\quad
\text{for all $ \iota \in \{1,...,d\}$}.
\end{equation}
We let
\begin{equation}\label{index_subsets_2}
I_{t_{1}} := 
I_{t}\setminus I_{t_{2}},
\qquad
I_{t_{2}} := 
\begin{cases}
\emptyset &\text{if $ \{x_{i} : i \in I_{t}\} \subseteq t_{1} $}\\
\{i \in I_{t} : x_{i} \in t_{2}\} &\text{otherwise}
\end{cases}
\end{equation}
and define $ \sons_{\ctI}(t) := \{t' \in \{t_{1}, t_{2}\} : I_{t'} \neq \emptyset\}. $
If $ \ell + 1 < \ell_{max}$ holds, we proceed with all sons $ t' \in \sons_{\ctI}(t) $ by recursion. Otherwise, we stop the recursion. 

By replacing $ I $ with $ J $ and $ (x_{i})_{i \in I} $ with $ (y_{j})_{j \in J} $ we analogously compute the tree $ \ctJ $ along with a corresponding hierarchical family $ (J_{s})_{s \in \ctJ} $ of non-empty subsets of $ J $.

 In conjunction with the families $ (I_{t})_{t\in \ctI} $ and $ (J_{s})_{s \in \ctJ} $, our trees $ \ctI $ and $ \ctJ $ form \emph{cluster trees} for $ I $ and $ J $, respectively, which implies
(cf. \cite[Corollary~3.9]{BO10})
\begin{align}\label{leafpart}
I = \dot{\bigcup_{t \in \Leaves_{\ctI}}} I_{t},
&&
J = \dot{\bigcup_{s \in \Leaves_{\ctJ}}} J_{s}.
\end{align}
Furthermore, by construction, we have
\begin{align}\label{leaves_lmax}
\Leaves_{\ctI} = \ctI^{(\ell_{max})},
&&
\Leaves_{\ctJ} = \ctJ^{(\ell_{max})}.
\end{align}

\begin{figure}[ht]
\begin{tikzpicture}[scale=0.7]
\begin{scope}[shift = {(3,0)}]
\draw[thick] (1,3) to[out=90, in = 70] (8,5);
\draw[thick] (8,5) to[out=250] (8,3);
\draw[thick] (8,3) to[out=-50,in=-90] (1,3);
\node[xshift=-0.3cm, yshift=-1.075cm] at (1,3) {\textcolor{red}{$ r^{(6)} $}};
\foreach \x in {1,2,3,4,5,6,7}{
	\draw[lightgray, thick, xshift = \x * 0.905cm] (1,1.18) rectangle (1.905, 1.82);
}
\foreach \x in {2,3,4,5,6,7}{
	\draw[lightgray, thick, xshift = \x * 0.905cm, yshift = 7 * 0.64cm] 
	(1,1.18) rectangle (1.905, 1.82);
}
\foreach \y in {1,2,3,4,5}{
	\draw[lightgray, thick, yshift = \y * 0.64cm] (1,1.18) rectangle (1.905, 1.82);
}
\foreach \y in {1,2,3,4,5,6}{
	\draw[lightgray, thick, xshift = 7 * 0.905cm, yshift = \y * 0.64cm] (1,1.18) rectangle (1.905, 1.82);
}

\node[xshift= 3 * 0.905cm + 0.3cm , yshift = 2 * 0.645cm] at (1,1.18){\textcolor{violet}{$ m_{t} $}};
\draw[blue, thick, xshift = 7 * 0.905cm, yshift = 5 * 0.64cm] (1,1.18) rectangle (1.905, 1.82);
\node[xshift= 7 * 0.905cm - 2.05cm, yshift = 5 * 0.645cm - 0.65cm] at (1,1.18) {\textcolor{blue}{$ t $}};
\node[xshift= 1 * 0.905cm + 0.05cm , yshift = 5 * 0.645cm - 0.325cm] at (1,1.18){\Large$ \Gamma $};

\draw[lightgray, thick, xshift = 1 * 0.905cm, yshift = 5 * 0.64cm] (1,1.18) rectangle (1.905, 1.82);
\draw[lightgray, thick, xshift = 1 * 0.905cm, yshift = 6 * 0.64cm] (1,1.18) rectangle (1.905, 1.82);
\draw[lightgray, thick, xshift = 2 * 0.905cm, yshift = 6 * 0.64cm] (1,1.18) rectangle (1.905, 1.82);
\draw[thick, xshift = 8.7cm] (1,3) to[out=90, in = 70] (8,5);
\draw[thick, xshift = 8.7cm] (8,5) to[out=250] (8,3);
\draw[thick, xshift = 8.7cm] (8,3) to[out=-50,in=-90] (1,3);
\node[xshift= 5.3cm, yshift=-1.1cm] at (1.5,3) {\textcolor{red}{$ B^{(6)} $}};
\foreach \x in {1,2,3,4,5,6,7}{
	\draw[lightgray, thick, xshift = \x * 0.905cm + 8.7cm] (1,1.18) rectangle (1.905, 1.82);
}
\foreach \x in {2,3,4,5,6,7}{
	\draw[lightgray, thick, xshift = \x * 0.905cm + 8.7cm, yshift = 7 * 0.64cm] 
	(1,1.18) rectangle (1.905, 1.82);
}
\foreach \y in {1,2,3,4,5}{
	\draw[lightgray, thick, xshift = 8.7cm, yshift = \y * 0.64cm] (1,1.18) rectangle (1.905, 1.82);
} 
\foreach \y in {1,2,3,4,5,6}{
	\draw[lightgray, thick, xshift = 7 * 0.905cm + 8.7cm, yshift = \y * 0.64cm] (1,1.18) rectangle (1.905, 1.82);
}
\node[xshift= 3 * 0.905cm + 6.4cm , yshift = 2 * 0.645cm] at (1,1.18){\textcolor{violet}{$ m_{t} $}};
\draw[blue, thick, xshift = 7 * 0.905cm + 8.7cm, yshift = 5 * 0.64cm] (0.9,1.05) rectangle (2.1,1.98);
\node[xshift= 7 * 0.905cm - 2.2cm + 6.1cm, yshift =  5 * 0.645cm - 0.6cm] at (1,1.18) {\textcolor{blue}{$ B_{t} $}};
\node[xshift= 1 * 0.905cm + 0.05cm + 6.1cm, yshift = 5 * 0.645cm - 0.325cm] at (1,1.18){\Large$ \Gamma $};
\draw[lightgray, thick, xshift = 1 * 0.905cm + 8.7cm, yshift = 5 * 0.64cm] (1,1.18) rectangle (1.905, 1.82);
\draw[lightgray, thick, xshift = 1 * 0.905cm + 8.7cm, yshift = 6 * 0.64cm] (1,1.18) rectangle (1.905, 1.82);
\draw[lightgray, thick, xshift = 2 * 0.905cm + 8.7cm, yshift = 6 * 0.64cm] (1,1.18) rectangle (1.905, 1.82);
\draw[lightgray,thick, xshift = 8.7cm] (1,1.18) rectangle (1.905,1.82);
\end{scope}
\draw[red, thick, xshift = 0cm, draw opacity = 0.0] (1,1.18) rectangle (1.905,1.82);
\begin{scope}[shift={(3,0)}]
\draw[red, thick, xshift = 0cm] (1,1.18) rectangle (1.905,1.82);
\draw[thick,-stealth, color=violet] (1.4525, 1.5) to (1.4525 + 7 * 0.905, 1.5 + 5 * 0.64);
\draw[red,thick, xshift = 8.7cm] (0.9,1.05) rectangle (2.1,1.98);
\draw[thick,-stealth, color=violet, xshift = 8.7cm] (1.4525, 1.5) to (1.4525 + 7 * 0.905, 1.5 + 5 * 0.64);
\end{scope}
\end{tikzpicture}
\caption{Relation between the boxes of $ \ctI $ and their support bounding
  boxes.}
\label{fi:supp_bboxes}
\end{figure}
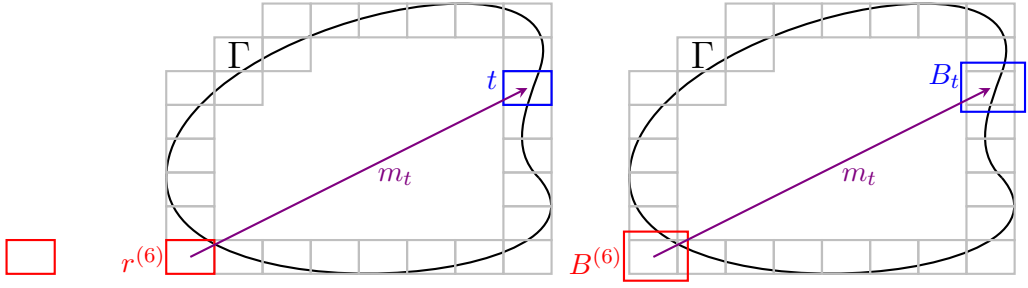

In the second phase we create corresponding \emph{support bounding boxes}
$ (B_{t})_{t\in\ctI} $ and $ (C_{s})_{s\in\ctJ} $ with (see
\cref{fi:supp_bboxes}) a property based on \cref{translational}: 
For every level $ \ell \in \{0,...,\ell_{max}\} $
we construct the smallest closed axis-parallel box $ B^{(\ell)} $ containing $ r^{(\ell)} $ and the shifted supports $ (\supp(\varphi_{i}) - m_{t})_{i\in I_{t}} $ for all $ t \in \ctI^{(\ell)} $. Likewise, for every level $ \ell \in \{0,...,\ell_{max}\} $
we construct the smallest closed axis-parallel box $ C^{(\ell)} $ containing $ r^{(\ell)} $ and the shifted supports $ (\supp(\psi_{j}) - m_{s})_{j\in J_{s}} $ for all $ s \in \ctJ^{(\ell)} $. According to (\ref{translational}), the properties
\begin{subequations}
\begin{equation}\label{support_bbox_I}
r^{(\ell)} \cup \Bigg(\bigcup_{i \in I_{t}}\supp(\varphi_{i}) - m_{t}\Bigg) \,\subseteq\, B^{(\ell)}
\qquad
\text{for all $t \in \ctI^{(\ell)} $},
\end{equation}
\begin{equation}\label{support_bbox_J}
r^{(\ell)}\cup \Bigg(\bigcup_{j \in J_{s}}\supp(\psi_{j}) - m_{s}\Bigg) \,\subseteq\, C^{(\ell)}
\qquad
\text{for all $s \in \ctJ^{(\ell)} $}
\end{equation}
\end{subequations}
then lead to
\begin{subequations}\label{shifted_support_bbox}
\begin{align}
t \,\cup\, \bigcup_{i \in I_{t}}\supp(\varphi_{i}) \,\subseteq\, B^{(\ell)} + m_{t} &\,=:\, B_{t} &\text{for all $t \in \ctI^{(\ell)} $},\\
s \,\cup\, \bigcup_{j \in J_{s}}\supp(\psi_{j}) \,\subseteq\, C^{(\ell)} + m_{s} &\,=:\, C_{s} &\text{for all $s \in \ctJ^{(\ell)} $}.
\end{align}
\end{subequations} 

%
%

In the third phase we construct an approximation of $ G $ by combining
the trees, boxes and index sets created before with tensor interpolation.
For this purpose we choose a degree $ \theta \in \bbbn $ and define
\begin{equation*} 
\Theta := \{0,...,\theta\}^{d}.
\end{equation*}
For every axis-parallel box $ Q \subset \bbbr^{d} $ let
\begin{align*}
  \interpol_{Q}\colon C(Q) &\to \Pi_\theta, &
     f &\mapsto \sum_{\nu\in\Theta} f(\xi_{Q,\nu}) \lagrange_{Q,\nu},
\end{align*}
be the corresponding tensor Chebyshev interpolation operator
mapping to the set $\Pi_\theta$ of tensor polynomials of degree
$\theta$ with associated interpolation points
$(\xi_{Q,\nu})_{\nu \in \Theta}$ and Lagrange polynomials
$(\lagrange_{Q,\nu})_{\nu \in \Theta}$.

We use the \emph{admissibility condition}
\begin{equation}\label{A_eta}
\begin{split}
A_{\eta} : T_{I} \times T_{J} &\rightarrow \{\true, \false\},\\
(t,s) &\mapsto
\begin{cases}
\true &\text{if $ \max\{\diam(B_{t}), \diam(C_{s})\} $}
\leq \eta\dist(B_{t}, C_{s})\\
\false &\text{if $ \max\{\diam(B_{t}), \diam(C_{s})\} $}
> \eta\dist(B_{t}, C_{s})             
\end{cases}
\end{split}
\end{equation}
with a given parameter $\eta\in\bbbr_{>0}$ to decide whether the
kernel function $g$ can be approximated in a domain $B_t\times C_s$.
For every pair $ (t,s) \in \ctI \times \ctJ$ with $ A_{\eta}(t,s) = \true $ the corresponding tensor interpolation error satisfies
\begin{equation*}
\|g - (\interpol_{B_{t}} \otimes \interpol_{C_{s}})[g]\|_{\infty, B_{t} \times C_{s}}
\leq \frac{\const_{in}}{\dist(B_{t}, C_{s})^{\sigma}}q_{\eta}^{\theta+1}
\end{equation*}
for a constant $ \const_{in} \in \bbbr_{>0} $ and $ q_{\eta} := \min\left\{\frac{\smallconst_{0}\eta}{\smallconst_{0}\eta+2}, \frac{\smallconst_{0}\eta}{4}\right\} $ (cf. \cite[Remark 4.23]{BO10}) and therefore
converges exponentially to zero with respect to the degree $ \theta $. The actual convergence rates that appear in practice are often considerably better than $ q_{\eta} $.

For every $ t \in \ctI $ let now
$ V_{t} \in \bbbr^{I_{t} \times \Theta}$ be given by
\begin{align*}
(V_{t})_{i\nu} &:= \int_{\Gamma}\varphi_{i}(x)\lagrange_{B_{t},\nu}(x) \, dx
& &\text{for all } i \in I_{t}, \nu \in \Theta,
\end{align*}
and for every $ s \in \ctJ $ let $ W_{s} \in \bbbr^{J_{s} \times \Theta} $
be given by
\begin{align*}
(W_{s})_{j\mu} &:= \int_{\Gamma}\psi_{j}(y)\lagrange_{C_{s},\mu}(y) \, dy
& &\text{for all } j \in J_{s}, \mu \in \Theta.
\end{align*}
Then for every pair $ (t,s) \in \ctI \times \ctJ$ with $ A_{\eta}(t,s) = \true $
we get, by using (\ref{eq:matrix_entries}) and (\ref{shifted_support_bbox}), for all $ i \in I_{t} , j \in J_{s}$
\begin{equation}\label{tensor_approx}
\begin{split}
G_{ij} 
&\approx 
\int_{\Gamma}\varphi_{i}(x)\int_{\Gamma}(\interpol_{B_{t}} \otimes \interpol_{C_{s}})[g](x,y)\psi_{j}(y)\, dy\, dx \\ 
&= \sum_{\nu \in \Theta}\sum_{\mu \in
  \Theta}\underbrace{\int_{\Gamma}\varphi_{i}(x)\lagrange_{B_{t},\nu}(x)\,
  dx}_{= (V_t)_{i\nu}} \, g(\xi_{B_{t},\nu},\xi_{C_{s},\mu}) \,
  \underbrace{\int_{\Gamma} \lagrange_{C_{s},\mu}(y)\psi_{j}(y) \,
  dy}_{= (W_s)_{j\mu}},
\end{split}
\end{equation}
which leads to
\begin{equation}\label{low-rank_fact}
G|_{I_{t} \times J_{s}} \approx V_{t}S_{(t,s)}W_{s}^{*} 
\end{equation}
with $ S_{(t,s)} \in \bbbr^{\Theta \times \Theta} $  given by
\begin{equation*}
(S_{(t,s)})_{\nu\mu} := g(\xi_{B_{t},\nu},\xi_{C_{s},\mu}) \qquad\textnormal{for all $ \nu \in \Theta, \mu \in \Theta $},
\end{equation*}
i.e., we have a factorized low-rank approximation of admissible
submatrices $G|_{I_t\times J_s}$.

We compute a tree $ \ctIJ \subseteq \ctI \times \ctJ$ defined by the following properties:
\begin{itemize}[leftmargin=0.4cm]
\item $ \Root_{\ctIJ} = (\Root_{\ctI},\Root_{\ctJ}) = (r^{(0)}, r^{(0)}) $;
\item For all $ (t,s) \in \ctIJ $ the sons are given by
  \begin{equation*}
    \sons_{\ctIJ}(t,s) = \begin{cases}
	\emptyset &\text{if } A_{\eta}(t,s) = \true,\\
	\sons_{\ctI}(t) \times \sons_{\ctJ}(s) &\text{otherwise.}
    \end{cases}
  \end{equation*}
  Note that this implies $\sons_{\ctIJ}(t,s)=\emptyset$ if
  $\sons_{\ctI}(t)=\emptyset$ or $\sons_{\ctJ}(s)=\emptyset$.
\end{itemize}
Its \emph{admissible leaves} and \emph{inadmissible leaves} are denoted by
\begin{align*}
\Leaves_{\ctIJ}^{+} &:= \{(t,s) \in \Leaves_{\ctIJ} : A_{\eta}(t,s) = \true\},\\
\Leaves_{\ctIJ}^{-} &:= \{(t,s) \in \Leaves_{\ctIJ} : A_{\eta}(t,s) = \false\}.
\end{align*}
Drawing on the properties of the tree $\ctI$ and $\ctJ$, we find
(cf. \cite[Corollary~3.15]{BO10})
\begin{equation*}
I \times J = \dot{\bigcup_{(t,s) \in \Leaves_{\ctIJ}}} (I_{t} \times J_{s})
= \dot{\bigcup_{(t,s) \in \Leaves^{+}_{\ctIJ}}} (I_{t} \times J_{s})
\quad \dot{\cup} \quad
\dot{\bigcup_{(t,s) \in \Leaves^{-}_{\ctIJ}}} (I_{t} \times J_{s}).
\end{equation*}
We can therefore use (\ref{low-rank_fact}) to construct an approximation $ \widetilde{G} \in \bbbr^{I \times J} $ of $ G $ by 
\begin{align*}
\widetilde{G}|_{I_{t} \times J_{s}} &:= 
V_{t}S_{(t,s)}W_{s}^{*} &\text{for all $ (t,s) \in \Leaves^{+}_{\ctIJ}$},\\ \nonumber
\widetilde{G}|_{I_{t} \times J_{s}} &:= G|_{I_{t} \times J_{s}} &\text{for all $ (t,s) \in \Leaves^{-}_{\ctIJ}. $}
\end{align*}
This means that for admissible blocks $(t,s)\in\Leaves^{+}_{\ctIJ}$,
only a small \emph{coupling matrix} $S_{(t,s)}\in\bbbr^{\Theta\times\Theta}$
has to be stored.

The \emph{row cluster basis} $(V_t)_{t\in\ctI}$ and the
\emph{column cluster basis} $(W_s)_{s\in\ctJ}$ can be stored efficiently
using the identity theorem for polynomials:
we have
\begin{equation*}
\lagrange_{B_{t},\nu} = \sum_{\mu \in \Theta}\lagrange_{B_{t},\nu}(\xi_{B_{t'},\mu})\lagrange_{B_{t'},\mu}
\qquad
\text{for all $ t \in \ctI\setminus\Leaves_{\ctI} , t' \in \sons_{\ctI}(t), \nu \in \Theta$}.
\end{equation*}
Consequently, for every $ t \in \ctI\setminus\Leaves_{\ctI}, t' \in \sons_{\ctI}(t), i \in I_{t'}, \nu \in \Theta $ we get
\begin{align*}
 (V_{t})_{i\nu} &= \int_{\Gamma}\varphi_{i}(x)\lagrange_{B_{t},\nu}(x) \, dx
				= \int_{\Gamma}\varphi_{i}(x)\sum_{\mu \in \Theta}\lagrange_{B_{t},\nu}(\xi_{B_{t'},\mu})\lagrange_{B_{t'},\mu}(x) \, dx\\
				&= \sum_{\mu \in \Theta}\int_{\Gamma}\varphi_{i}(x)\lagrange_{B_{t'},\mu}(x) \, dx \, \, \lagrange_{B_{t},\nu}(\xi_{B_{t'},\mu}) 
				= \sum_{\mu \in \Theta}(V_{t'})_{i\mu} \, \, \lagrange_{B_{t},\nu}(\xi_{B_{t'},\mu})
\end{align*}
and therefore
\begin{equation}\label{transfer_row}
V_{t}|_{I_{t'} \times \Theta} = V_{t'}E_{t'}
\end{equation}
with $ E_{t'} \in \bbbr^{\Theta \times \Theta} $ given by
\begin{equation*}
(E_{t'})_{\mu\nu} := \lagrange_{B_{t},\nu}(\xi_{B_{t'},\mu})
\qquad 
\text{for all $ \mu \in \Theta, \nu \in \Theta$}.
\end{equation*}
Since we also have
$ I_{t} = \dot{\bigcup}_{t' \in \sons_{\ctI}(t)}\,I_{t'}  $ 
for all $ t \in \ctI\setminus\Leaves_{\ctI} $ by construction,
we do not have to store all the matrices $ (V_{t})_{t \in \ctI} $
explicitly.
It is sufficient to store only the
\emph{leaf matrices} $ (V_{t})_{t \in \Leaves_{\ctI}} $ and the
\emph{transfer matrices} $ (E_{t})_{t \in \ctI\setminus\{\Root_{\ctI}\}} $ instead.

So far, we have not taken advantage of the translation-invariance
\cref{eq:translation} of the kernel function.
It allows us to avoid storing all the coupling matrices
$(S_{b})_{b\in \Leaves_{\ctIJ}^{+}} $ individually: 
\cref{shifted_support_bbox} implies
\begin{subequations}\label{shifted interpoints}
\begin{align}
\label{shifted interpoints B_t}
\xi_{B_{t},\nu} &= \xi_{B^{(\ell)},\nu} + m_{t}  &\text{for all $ t \in \ctI^{(\ell)}, \nu \in \Theta $},\\
\label{shifted interpoints C_s}
\xi_{C_{s},\mu} &= \xi_{C^{(\ell)},\mu} + m_{s}  &\text{for all $ s \in \ctJ^{(\ell)}, \mu \in \Theta $}
\end{align}
\end{subequations}
on every level $ \ell \in \{0,...,\ell_{max}\} $ while, 
due to \cref{leaves_lmax}, the tree $ \ctIJ $ satisfies
\begin{equation}
\label{all same level}
\level_{\ctIJ}(t,s) = \level_{\ctI}(t) = \level_{\ctJ}(s) 
\qquad \quad
\text{for all $ (t,s) \in \ctIJ $}.
\end{equation}
\noindent By combining \cref{all same level} with \cref{shifted interpoints} and
\cref{eq:translation} we get
\begin{align*}
(S_{(t,s)})_{\nu\mu} 
&= g(\xi_{B_{t}, \nu}, \xi_{C_{s},\mu}) 
= g(\xi_{B^{(\ell)}, \nu} + m_{t}, \xi_{C^{(\ell)},\mu} + m_{s}) \\
&= g(\xi_{B^{(\ell)}, \nu} + m_{t} - m_{t}, \xi_{C^{(\ell)},\mu} + m_{s} - m_{t})\\
&= g(\xi_{B^{(\ell)}, \nu}, \xi_{C^{(\ell)},\mu} - (m_{t} - m_{s}))
\qquad
\text{for all $ (t,s) \in \Leaves_{\ctIJ}^{+} \cap \ctIJ^{(\ell)} $}, \nu, \mu \in \Theta
\end{align*}
on every level $ \ell \in \{0,...,\ell_{max}\} $.
Hence, on every level $ \ell \in \{0,...,\ell_{max}\} $ the implication
\begin{equation}\label{tinv_coupling}
m_{t} - m_{s} = m_{\tilde{t}} - m_{\tilde{s}} 
\enspace \implies \enspace
S_{(t,s)} = S_{(\tilde{t},\tilde{s})}
\end{equation}
holds for all $ (t,s), (\tilde{t},\tilde{s}) \in \Leaves_{\ctIJ}^{+}\cap\ctIJ^{(\ell)} $.

Moreover, due to our regular construction of the cluster tree, we also do
not have to store all the transfer matrices individually: 
On every level $ \ell \in \{0,...,\ell_{max}-1\} $ we have,
according to \cref{shifted interpoints B_t},
\begin{align*}
\lagrange_{B_{t},\nu}(x) &= \lagrange_{B^{(\ell)},\nu}(x - m_{t}) &
  &\text{ for all } t \in \ctI^{(\ell)}, \nu \in \Theta, x \in \bbbr^{d},
\end{align*}
which implies
\begin{align*}
(E_{t'})_{\mu\nu} 
&= \lagrange_{B_{t},\nu}(\xi_{B_{t'},\mu}) 
= \lagrange_{B^{(\ell)},\nu}(\xi_{B_{t'},\mu} - m_{t}) \\
&= \lagrange_{B^{(\ell)},\nu}(\xi_{B^{(\ell+1)},\mu} + m_{t'} - m_{t})\\
&= \lagrange_{B^{(\ell)},\nu}(\xi_{B^{(\ell+1)},\mu} - (m_{t} - m_{t'}))
  \qquad \text{ for all } t \in \ctI^{(\ell)}, t' \in \sons_{\ctI}(t), \nu, \mu \in \Theta.
\end{align*}
Hence, on every level $\ell \in \{0,...,\ell_{max}-1\}$, we have
\begin{equation}\label{tinv_row}
m_{t} - m_{t'} = m_{\tilde{t}} - m_{\tilde{t}'} \enspace \implies \enspace
E_{t'} = E_{\tilde{t}'} 
\end{equation}
for all $ t,\tilde{t}\in\ctI^{(\ell)},
t'\in\sons_{\ctI}(t),\ \tilde{t}'\in\sons_{\ctI}(\tilde{t}) $.
The column cluster basis $(W_s)_{s\in\ctJ}$ and the corresponding
transfer matrices $(F_s)_{s\in\ctJ}$ have similar properties.

The implications \cref{tinv_coupling} and \cref{tinv_row} allow us
to significantly reduce the storage requirements of the $\htwo$-matrix
representation by taking advantage of the translation-invariance
of the kernel function.

\section{Complexity estimates}
\label{se:complexity}

In the following, let 
\begin{align}\label{nk}
  n := \max\{|I|,|J|\},
  &&
  k := |\Theta| = (\theta+1)^{d}.
\end{align}
For every $ \iota \in \{1,...,d\} $ let $ e_{\iota} \in \bbbr^{d}$ denote the $ \iota$-th canonical unit vector.

%
%
\begin{lemma}[Storage requirements of the leaf matrices]
\label{leaf storage}
Storing the leaf matrices $(V_{t})_{t \in \Leaves_{\ctI}}$ and
$(W_{s})_{s \in \Leaves_{\ctJ}}$ requires not more than $ 2kn $ units
of storage.
\end{lemma}
\textit{Proof.} Combining \cref{nk} with \cref{leafpart} implies that storing the matrices $ (V_{t})_{t \in \Leaves_{\ctI}} $ and $ (W_{s})_{s \in \Leaves_{\ctJ}} $ directly requires not more than
\begin{equation*}
\sum_{t \in \Leaves_{\ctI}} |I_{t}| \, |\Theta|
\enspace +  \sum_{s \in \Leaves_{\ctJ}} |J_{s}| \, |\Theta|
=
k\bigg(\sum_{t \in \Leaves_{\ctI}}|I_{t}| \enspace + \sum_{s \in \Leaves_{\ctJ}}|J_{s}|\bigg)
= k(|I| + |J|) \leq 2kn
\end{equation*}
units of storage. $\hfill \Box  $

\begin{figure}[ht]
\begin{center}
\begin{tikzpicture}[scale = 0.7]
\draw[thick] (1,3) to[out=90, in = 70] (8,5);
\draw[thick] (8,5) to[out=250] (8,3);
\draw[thick] (8,3) to[out=-50,in=-90] (1,3);
\node[xshift=-0.35cm, yshift=-0.85cm] at (2,2.7) {\textcolor{red}{$ r^{(5)} $}};
\foreach \x in {1,2,3,4,5,6,7}{
	\draw[lightgray, thick, xshift= 1.0cm + \x * 0.905cm, yshift = 1.18cm] 
	(0,0) rectangle (0.905,1.28);
}
\foreach \x in {1,2,3,4,5,6,7}{
	\draw[lightgray, thick, xshift= 1.0cm + \x * 0.905cm, yshift = 1.18cm + 3 * 1.28 cm] 
	(0,0) rectangle (0.905,1.28);
}
\foreach \y in {1,2,3}{
	\draw[lightgray, thick, xshift= 1.0cm + 7 * 0.905cm, yshift = 1.18cm + \y * 1.28 cm] 
	(0,0) rectangle (0.905,1.28);
}
\foreach \y in {1,2}{
	\draw[lightgray, thick, xshift= 1.0cm, yshift = 1.18cm + \y * 1.28 cm] 
	(0,0) rectangle (0.905,1.28);
}
\draw[lightgray, thick, xshift= 1.0cm + 0.905cm, yshift = 1.18cm + 2 * 1.28 cm] 
(0,0) rectangle (0.905,1.28);
\draw[red, thick] (1,1.18) rectangle (1.905,2.46);
\draw[red, thick, xshift= 1.0cm + 1 * 0.905cm, yshift = 1.18cm + 0 * 1.28 cm] (0,0) rectangle (0.905,1.28);
\draw[blue, thick, xshift= 1.0cm + 2 * 0.905cm, yshift = 1.18cm + 3 * 1.28 cm] (0,0) rectangle (0.905,1.28);
\draw[blue, thick, xshift= 1.0cm + 3 * 0.905cm, yshift = 1.18cm + 3 * 1.28 cm] (0,0) rectangle (0.905,1.28);
\node[xshift=0.0cm, yshift=-1.35cm] at (2,2.7) {\Large\textcolor{red}{$ r^{(4)} $}};
\node[xshift= 1 * 0.905cm + 0.05cm , yshift = 5 * 0.645cm - 0.325cm] at (1,1.18){\Large$ \Gamma $};
\node[xshift= 3 * 0.905cm - 0.8cm , yshift = 6 * 0.645cm -0.075cm] at (1,1.18)
{\textcolor{blue}{$ t $}};
\node[xshift= 2 * 0.905cm - 0.45cm , yshift = 2 * 0.645cm - 0.1cm] at (1,1.18)
{\textcolor{violet}{$ m_{t} $}};
\node[xshift= 2 * 0.905cm - 1cm , yshift = 4 * 0.645cm - 0.4cm] at (1,1.18)
{\textcolor{teal}{$ m_{t_{1}} $}};
\node[xshift= 4 * 0.905cm - 1.45cm , yshift = 4 * 0.645cm - 0.2cm] at (1,1.18)
{\textcolor{teal}{$ m_{t_{2}} $}};
\draw[-stealth, teal, thick] (1.4525,1.82) to (1.4525 + 2 * 0.905, 1.82 + 3 * 1.28);
\draw[-stealth, teal, thick] (1.4525,1.82) to (1.4525 + 3 * 0.905, 1.82 + 3 *1.28);
\draw[-stealth, violet, thick] (1.905, 1.82) to (1.905 + 2 * 0.905, 1.82 + 3 * 1.28);
\end{tikzpicture}
\begin{tikzpicture}[scale = 0.7]
\draw[thick] (1,3) to[out=90, in = 70] (8,5);
\draw[thick] (8,5) to[out=250] (8,3);
\draw[thick] (8,3) to[out=-50,in=-90] (1,3);
\node[xshift=-0.35cm, yshift=-0.85cm] at (2,2.7) {\textcolor{red}{$ r^{(5)} $}};
\foreach \x in {1,2,3,4,5,6,7}{
	\draw[lightgray, thick, xshift= 1.0cm + \x * 0.905cm, yshift = 1.18cm] 
	(0,0) rectangle (0.905,1.28);
}
\foreach \x in {1,2,3,4,5,6,7}{
	\draw[lightgray, thick, xshift= 1.0cm + \x * 0.905cm, yshift = 1.18cm + 3 * 1.28 cm] 
	(0,0) rectangle (0.905,1.28);
}
\foreach \y in {1,2,3}{
	\draw[lightgray, thick, xshift= 1.0cm + 7 * 0.905cm, yshift = 1.18cm + \y * 1.28 cm] 
	(0,0) rectangle (0.905,1.28);
}
\foreach \y in {1,2}{
	\draw[lightgray, thick, xshift= 1.0cm, yshift = 1.18cm + \y * 1.28 cm] 
	(0,0) rectangle (0.905,1.28);
}
\draw[lightgray, thick, xshift= 1.0cm + 0.905cm, yshift = 1.18cm + 2 * 1.28 cm] 
(0,0) rectangle (0.905,1.28);
\draw[red, thick] (1,1.18) rectangle (1.905,2.46);
\draw[red, thick, xshift= 1.0cm + 1 * 0.905cm, yshift = 1.18cm + 0 * 1.28 cm] (0,0) rectangle (0.905,1.28);
\draw[blue, thick, xshift= 1.0cm + 4 * 0.905cm, yshift = 1.18cm + 3 * 1.28 cm] (0,0) rectangle (0.905,1.28);
\draw[blue, thick, xshift= 1.0cm + 5 * 0.905cm, yshift = 1.18cm + 3 * 1.28 cm] (0,0) rectangle (0.905,1.28);
\node[xshift=0.0cm, yshift=-1.35cm] at (2,2.7) {\Large\textcolor{red}{$ r^{(4)} $}};
\node[xshift= 1 * 0.905cm + 0.05cm , yshift = 5 * 0.645cm - 0.325cm] at (1,1.18){\Large$ \Gamma $};
\node[xshift= 2 * 0.905cm + 0.0cm , yshift = 2 * 0.645cm - 0.1cm] at (1,1.18)
{\textcolor{violet}{$ m_{\tilde{t}} $}};
\node[xshift= 3 * 0.905cm - 1cm , yshift = 4 * 0.645cm - 0.15cm] at (1,1.18)
{\textcolor{teal}{$ m_{\tilde{t}_{1}} $}};
\node[xshift= 4 * 0.905cm - 0.5cm , yshift = 4 * 0.645cm - 0.25cm] at (1,1.18)
{\textcolor{teal}{$ m_{\tilde{t}_{2}} $}};
\node[xshift= 5 * 0.905cm - 1.35cm , yshift = 6 * 0.645cm -0.05cm] at (1,1.18)
{\textcolor{blue}{$ \tilde{t} $}};
\draw[-stealth, teal, thick] (1.4525,1.82) to (1.4525 + 4 * 0.905, 1.82 + 3 * 1.28);
\draw[-stealth, teal, thick] (1.4525,1.82) to (1.4525 + 5 * 0.905, 1.82 + 3 *1.28);
\draw[-stealth, violet, thick] (1.905, 1.82) to (1.905 + 4 * 0.905, 1.82 + 3 * 1.28);
\end{tikzpicture}
\end{center}
\caption{The translations of the boxes of $ \ctI $ and the translations of their sons.}
\label{sondiffvecfig}
\end{figure}
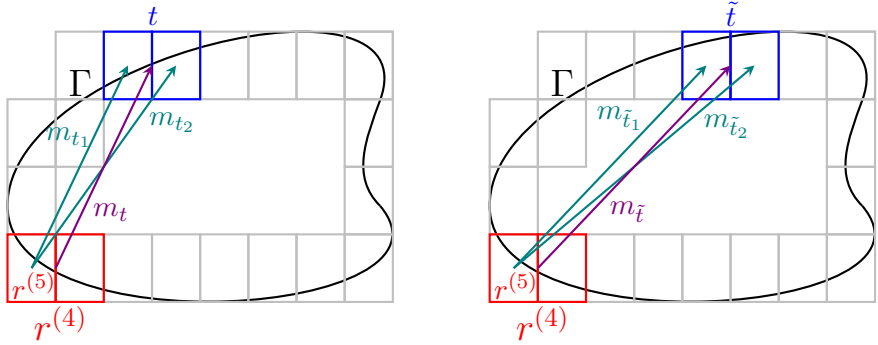

%
%
\begin{lemma}\label{sondiffvec}
  On every level $ \ell \in \{0, ..., \ell_{max} - 1\} $ we have
  \begin{align*}
  \{m_{t} - m_{t'} : t \in \ctI^{(\ell)}, t' \in \sons_{\ctI}(t)\}
  &\subseteq
  \{0, -\delta^{(\ell+1)}_{\iota^{(\ell)}}e_{\iota^{(\ell)}}\}, 
   \\
  \{m_{s} - m_{s'} : s \in \ctJ^{(\ell)}, s' \in \sons_{\ctJ}(s)\}
  &\subseteq \{0, -\delta^{(\ell+1)}_{\iota^{(\ell)}}e_{\iota^{(\ell)}}\}.
  \end{align*}
\end{lemma}
\textit{Proof.} 
Let $ \ell \in \{0,...,\ell_{max}-1\} $, let $ t \in \ctI^{(\ell)} $ and
let $ t' \in \sons_{\ctI}(t) $.
Then $ p_{t'} \in \bbbn_{0}^{d}$ satisfies $ p_{t',\iota} = p_{t,\iota} $
for all $ \iota \in \{1,...,d\}\setminus\{\iota^{(\ell)}\}$ and
$ p_{t',\iota^{(\ell)}} \in \{2p_{t,\iota^{(\ell)}}, 2p_{t,\iota^{(\ell)}} +
1\} $,
according to \cref{displacement_sons}.
Due to \cref{translational} and \cref{halved_ref_box} we therefore get, 
as illustrated in Figure \ref{sondiffvecfig}, 
$(m_{t} - m_{t'})_{\iota} = 0$ for all
$\iota \in \{1,...,d\}\setminus\{\iota^{(\ell)}\}$ and
\begin{align*}
(m_{t} - m_{t'})_{\iota^{(\ell)}} 
&= p_{t,\iota^{(\ell)}}\delta^{(\ell)}_{\iota^{(\ell)}} - p_{t',\iota^{(\ell)}}\delta^{(\ell+1)}_{\iota^{(\ell)}}\\
&\in 
\{p_{t,\iota^{(\ell)}}\delta^{(\ell)}_{\iota^{(\ell)}} - 2p_{t,\iota^{(\ell)}}\delta^{(\ell+1)}_{\iota^{(\ell)}}, \, p_{t,\iota^{(\ell)}}\delta^{(\ell)}_{\iota^{(\ell)}} - 2p_{t,\iota^{(\ell)}}\delta^{(\ell+1)}_{\iota^{(\ell)}}
- \delta^{(\ell+1)}_{\iota^{(\ell)}}\}\\
&=
\{p_{t,\iota^{(\ell)}}\delta^{(\ell)}_{\iota^{(\ell)}} - p_{t,\iota^{(\ell)}}\delta^{(\ell)}_{\iota^{(\ell)}},  \,  p_{t,\iota^{(\ell)}}\delta^{(\ell)}_{\iota^{(\ell)}} - p_{t,\iota^{(\ell)}}\delta^{(\ell)}_{\iota^{(\ell)}}
- \delta^{(\ell+1)}_{\iota^{(\ell)}}\}\\
&= \{0, - \delta^{(\ell+1)}_{\iota^{(\ell)}}\},
\end{align*} 
which proves the first statement.  The second statement can be proven similarly.
 $\hfill \Box $ 

%
%

%
%
\begin{theorem}[Storage requirements of the transfer matrices]
\label{transfer storage}
The transfer matrices $(E_{t})_{t\in\ctI\setminus\{\Root_{\ctI}\}}$
and $(F_{s})_{s\in\ctJ\setminus\{\Root_{\ctJ}\}}$ require not more than 
$4k^{2}\ell_{max}$ units of storage.
\end{theorem}
\textit{Proof.} Combining the implication \cref{tinv_row} 
with \cref{sondiffvec} and \cref{nk} directly implies that for
representing $(E_{t})_{t\in\ctI\setminus\{\Root_{\ctI}\}} $ and
$ (F_{s})_{s\in\ctJ\setminus\{\Root_{\ctJ}\}}$ not more than
\begin{align*}
&\,\quad|\Theta|^{2}\sum_{\ell = 0}^{\ell_{max}-1}\Big|\{m_{t} - m_{t'} : t \in \ctI^{(\ell)}, t' \in \sons_{\ctI}(t)\}\Big|\\
&+|\Theta|^{2}\sum_{\ell = 0}^{\ell_{max}-1}\Big|\{m_{s} - m_{s'} : s \in \ctJ^{(\ell)}, s' \in \sons_{\ctJ}(s)\}\Big|\\
&\leq |\Theta|^{2}\sum_{\ell = 0}^{\ell_{max}-1}4
\,=\, 4k^{2}\ell_{max}
\end{align*} 
units of storage are needed. $ \hfill \Box $

In order to obtain further results, we now consider not only
the boxes appearing as nodes in our trees $ \ctI $ and $ \ctJ $, but
the infinitely many boxes belonging to
$\that := \bigcup_{\ell=0}^{\ell_{max}}\widehat{T}^{(\ell)}$ given
by (cf. \cref{translational})
\begin{align}\label{all_placements}
\widehat{T}^{(\ell)} := \{r^{(\ell)} + \delta^{(\ell)} \odot p : p \in \bbbz^{d}\}
\qquad \text{for all $ \ell \in \{0,...,\ell_{max}\} $}.
\end{align}
For all $ \ell \in \{0,...,\ell_{max}\} $ let $ (p_{t})_{t\in\widehat{T}^{(\ell)}} $ be defined by
\begin{align}\label{general_placement}
p_{t} &\in \bbbz^{d}, &
t &=
r^{(\ell)} + \underbrace{\delta^{(\ell)} \odot p_{t}}_{=: m_{t}} 
= r^{(\ell)} + m_{t} &
&\text{ for all } t \in \widehat{T}^{(\ell)}
\end{align}
and let (cf. \cref{shifted_support_bbox})
\begin{subequations}\label{shifted_ref_supp_bbox}
\begin{align}
\label{shifted_ref_supp_bbox_B_t}
B_{t} &:= B^{(\ell)} + \delta^{(\ell)} \odot p_{t}
= B^{(\ell)} + m_{t} &
  &\text{ for all } t \in \that^{(\ell)},\\
\label{shifted_ref_supp_bbox_C_s}
C_{s} &:= C^{(\ell)} + \delta^{(\ell)} \odot p_{s}
= C^{(\ell)} + m_{s} &
  &\text{ for all } s \in \that^{(\ell)}.
\end{align}
\end{subequations}
Based on (\ref{A_eta}), let
\begin{equation}\label{Ahat_eta}
\begin{split}
\ahat_{\eta} : \that \times \that &\rightarrow \{\true, \false\}, \\
(t,s) &\mapsto
\begin{cases}
\true &\text{if $ \max\{\diam(B_{t}), \diam(C_{s})\} $}
\leq \eta\dist(B_{t}, C_{s})\\
\false &\text{if $ \max\{\diam(B_{t}), \diam(C_{s})\} $}
> \eta\dist(B_{t}, C_{s})             
\end{cases}.
\end{split}
\end{equation}

We base the following analysis on one key assumption.

%
%
\begin{assumption}[Diameters of supports]
We assume that there is a constant $\const_{bb} \in \bbbr_{\geq1}$
satisfying
\begin{align}\label{diameter_outer_inner}
  \diam(B^{(\ell)}) &\leq \const_{bb}\diam(r^{(\ell)}) \qquad\text{ and} &
  \diam(C^{(\ell)}) &\leq \const_{bb}\diam(r^{(\ell)})
\end{align}
for all $ \ell \in \{0,...,\ell_{max}\}$.
\end{assumption}

We will see later on that under normal circumstances, due to our
choice of $\ell_{max}$, this condition is actually fulfilled
since we do not allow the boxes to become too small
(see Remark \ref{assumption_true}).

Let $\lambda_{d}$ denote the $d$-dimensional Lebesgue measure. 

%
%
\begin{lemma}
There is a constant $\const_{dv} \in \bbbr_{> 1}$ such
that
\begin{align}\label{diameter_volume}
  \diam(r^{(\ell)})^{d} &\leq \const_{dv}\lambda_{d}(r^{(\ell)}) &
  &\text{ for all } \ell \in \{0,...,\ell_{max}\}.
\end{align}
\end{lemma}
\textit{Proof.} According to \cref{regular_splitting} and \cref{halved_ref_box}, we have
\begin{equation*}
\diam(r^{(\ell)}) = \tfrac{1}{2}\diam(r^{(\ell-d)}) \quad \textnormal{and} \quad
\lambda_{d}(r^{(\ell)}) = \tfrac{1}{2^{d}}\lambda_{d}(r^{(\ell-d)})
\qquad \textnormal{for all $ \ell \in \{d,...,\ell_{max}\} $}
\end{equation*}
and therefore
\begin{equation*}
\frac{\diam(r^{(\ell)})^{d}}{\lambda_{d}(r^{(\ell)})}
= \frac{\diam(r^{(\ell-d)})^{d}}{\lambda_{d}(r^{(\ell-d)})}
\qquad \textnormal{for all $ \ell \in \{d,...,\ell_{max}\} $}.
\end{equation*}
This implies
\begin{equation*}
\left\{\frac{\diam(r^{(\ell)})^{d}}{\lambda_{d}(r^{(\ell)})} : \ell \in \{0,...,\ell_{max}\}\right\}
= \left\{\frac{\diam(r^{(\ell)})^{d}}{\lambda_{d}(r^{(\ell)})} : \ell \in \{0,...,d-1\}\right\} =: R.
\end{equation*}
Hence, $ \const_{dv} := \max R$ satisfies \cref{diameter_volume}. $ \hfill \Box $

\begin{figure}[ht]\label{inadm_cont_ball}
	\begin{center}
		\begin{tikzpicture}[scale = 0.7, baseline={([yshift=-.8ex]current bounding box.center)}]
		\draw[lightgray, thick] (1,1.18) rectangle (1.905,1.82);
		\node[xshift = -0.03cm, yshift = 0.93cm] at (1.81,0) {\textcolor{gray}{\tiny\thicklines$ t $}};
		\foreach \x in {-1, 0, 1} {
			\foreach \y in {-2, -1, 0, 1, 2}{
				\draw[lightgray, thick, xshift = \x * 0.905cm, yshift = \y * 0.64cm] (1,1.18) rectangle (1.905,1.82);
			}
		}
		\foreach \x in {-2, 2} {
			\foreach \y in {-2, -1, 0, 1, 2}{
				\draw[lightgray, thick, xshift = \x * 0.905cm, yshift = \y * 0.64cm] (1,1.18) rectangle (1.905,1.82);
			}
		}
		\foreach \x in {-1, 0, 1} {
			\foreach \y in {-3, 3}{
				\draw[lightgray, thick, xshift = \x * 0.905cm, yshift = \y * 0.64cm] (1,1.18) rectangle (1.905,1.82);
			}
		}
		\foreach \y in {-4, 4}{
			\draw[lightgray, thick, xshift = 0 * 0.905cm, yshift = \y * 0.64cm] (1,1.18) rectangle (1.905,1.82);
		}
		\draw[blue, thick] (0.9,1.05) rectangle (2.1,1.98);
		\node[xshift = 0.25cm, yshift = 0.7cm] at (2.0,-0.18) {\small\textcolor{blue}{$ B_{t} $}};
		\node[yshift=-0.015cm] at (1.5, 1.515) {\textcolor{black}{\LARGE$\cdot$}}; 
		\node[xshift=0.33cm, yshift=0.53cm] at (0.9, 0.645) {\tiny$\beta_{t}$};
		\draw[thick] (1.5, 1.515) to (4.54,1.515);
		\draw[thick] (1.5, 1.515) circle (3.04);
		\node[xshift=1.7cm, yshift=0.84cm] at (0.9, 0.645) {\small$\rho^{(6)}$};
		\node[xshift=-2.1cm, yshift=-1cm] at (0.9, 0.645) {$ \Bsf_{2}[\beta_{t}, \rho^{(6)}] $};
		\end{tikzpicture}
	\end{center}
	\caption{For $ t \in \widehat{T}^{(6)} $ the ball $
          \Bsf_{2}[\beta_{t}, \rho^{(6)}] $ contains all inadmissible $ s \in \widehat{T}^{(6)} $.}
\end{figure}
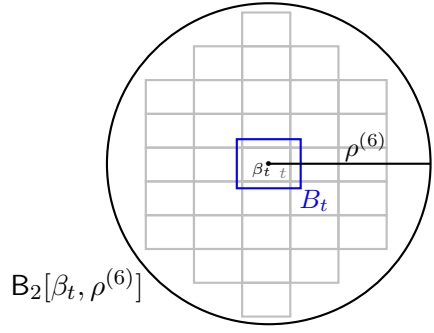

%
%
\begin{lemma}[Sparsity]\label{general sparsity}
We define
\begin{equation*}
 \const_{sp} := 2^{-d}\left(3 + 2\eta^{-1}\right)^{d}
               \const_{bb}^{d}\omega_{d}\const_{dv},
\end{equation*}
where $ \omega_{d} := \lambda_{d}(\{z \in \bbbr^{d} : \|z\|_{2} \leq 1\})$
is the volume of the $ d $-dimensional unit ball.

For every $ \ell \in \{0,...,\ell_{max}\} $ we have
\begin{align*}
  \left|\{s \in \widehat{T}^{(\ell)} : \ahat_{\eta}(t,s) = \false\}\right| 
  &\leq \const_{sp} &
  &\text{ for all } t \in \widehat{T}^{(\ell)},\\
  \left|\{t \in \widehat{T}^{(\ell)} : \ahat_{\eta}(t,s) = \false\}\right|
  &\leq \const_{sp} &
  &\text{ for all } s \in \widehat{T}^{(\ell)}.
\end{align*}
\end{lemma}
\textit{Proof.}
Let $ \ell \in \{0,...,\ell_{max}\} $, let $ t \in \widehat{T}^{(\ell)} $,
let $H_{t} := \{s \in \widehat{T}^{(\ell)} : \ahat_{\eta}(t,s) = \false\}$
and let $ \beta_{t} $ denote the midpoint of $ B_{t} $.
For every $ s \in H_{t} $ we have according to \cref{Ahat_eta}
\begin{equation*}
  \dist(B_{t}, C_{s}) < \eta^{-1}\max\{\diam(B_{t}), \diam(C_{s})\}.
\end{equation*} 
Thus for every $ s \in H_{t} $ there is at least one pair $ (x_{s}, y_{s}) \in B_{t} \times C_{s} $ satisfying
$ \|x_{s} - y_{s}\|_{2} < \eta^{-1}\max\{\diam(B_{t}),\diam(C_{s})\} $,
which, according to (\ref{shifted_ref_supp_bbox}), implies
\begin{align*}
\|z - \beta_{t}\|_{2} &\leq \|z - y_{s}\|_{2} + \|y_{s} - x_{s}\|_{2} + \|x_{s} - \beta_{t}\|_{2}\\
&< \diam(C_{s}) + \eta^{-1}\max\{\diam(B_{t}),\diam(C_{s})\} + \tfrac{1}{2}\diam(B_{t})\\
&\leq \tfrac{1}{2}(3 + 2\eta^{-1})\max\{\diam(B_{t}),\diam(C_{s})\}\\
&= \tfrac{1}{2}(3 + 2\eta^{-1})\max\{\diam(B^{(\ell)}),\diam(C^{(\ell)})\}
=: \rho^{(\ell)} \qquad \text{for all $ z \in C_{s} $},
\end{align*}
which, combined with
\cref{support_bbox_J}, \cref{general_placement} and
\cref{shifted_ref_supp_bbox_C_s}, implies (cf. Figure \ref{inadm_cont_ball})
\begin{equation}\label{inner_contained}
s \subseteq C_{s} \,\subset\, \Bsf_{2}[\beta_{t}, \rho^{(\ell)}] := \{z \in \bbbr^{d} : \|z - \beta_{t}\|_{2} \leq \rho^{(\ell)}\}.
\end{equation}
Moreover, according to \cref{all_placements}, we have (cf. Figure \ref{inadm_cont_ball})
\begin{equation}\label{gt1_null_set}
\lambda_{d}\big(\big\{z \in \Bsf_{2}[\beta_{t}, \rho^{(\ell)}] : |\{s \in H_{t} : z \in s\}| > 1\big\}\big) = 0.
\end{equation}
By using \cref{diameter_volume}, \cref{all_placements},
\cref{inner_contained}, \cref{gt1_null_set} and \cref{diameter_outer_inner}
we now get
\begin{align*}
\diam(r^{(\ell)})^{d}|H_{t}| 
&\leq \const_{dv}\lambda_{d}(r^{(\ell)})|H_{t}| 
= \const_{dv}\sum_{s \in H_{t}}\lambda_{d}(s) 
= \const_{dv}\sum_{s \in H_{t}}\int_{\Bsf_{2}[\beta_{t}, \rho^{(\ell)}]}\mathbbm{1}_{s}(z) \,dz \\
&= \const_{dv}\int_{\Bsf_{2}[\beta_{t}, \rho^{(\ell)}]}\sum_{s \in H_{t}}\mathbbm{1}_{s}(z) \, \, dz
\leq \const_{dv}\int_{\Bsf_{2}[\beta_{t}, \rho^{(\ell)}]}\mathbbm{1}_{\Bsf_{2}[\beta_{t}, \rho^{(\ell)}]}(z) \,dz \\
&= \const_{dv}\lambda_{d}(\Bsf_{2}[\beta_{t}, \rho^{(\ell)}])
= \const_{dv}\omega_{d}(\rho^{(\ell)})^{d} \\
&= \const_{dv}\omega_{d}2^{-d}(3 + 2\eta^{-1})^{d}\max\{\diam(B^{(\ell)}),\diam(C^{(\ell)})\}^{d}\\
&\leq \const_{dv}\omega_{d}2^{-d}(3 + 2\eta^{-1})^{d}\const_{bb}^{d}\diam(r^{(\ell)})^{d},
\end{align*}
which proves the first statement. The second statement can be proven similarly.
$ \hfill \Box $ 

%
%

%
%
\begin{lemma}[Difference bound for $ \ctIJ $]\label{diff_bound}
  On every level $ \ell \in \{0,...,\ell_{max}\} $ we have
  \begin{equation*}
  \left|\{m_{t} - m_{s} : (t,s) \in \ctIJ^{(\ell)}\}\right| \leq 3\const_{sp}.
  \end{equation*}	
\end{lemma}
\textit{Proof.} By definition, we have
$ \ctIJ^{(0)} = \{(\Root_{\ctI},\Root_{\ctJ})\} = \{(r^{(0)}, r^{(0)})\}$
and therefore
\begin{equation*}
  \{m_{t} - m_{s} : (t,s) \in \ctIJ^{(0)}\}
  = \{m_{r^{(0)}} - m_{r^{(0)}}\} = \{0\}.
\end{equation*}
Let $ \ell \in \{1,...,\ell_{max}\} $ and let $ (t,s) \in \ctIJ^{(\ell)} $.
Then, by the definition of $ \ctIJ $ and due to \cref{all same level},
there is a pair
$(\check{t},\check{s}) \in \ctIJ^{(\ell-1)} \subseteq \ctI^{(\ell-1)}
\times \ctJ^{(\ell-1)}$ satisfying $t \in \sons_{\ctI}(\check{t})$,
$s \in \sons_{\ctJ}(\check{s})$ and $A_{\eta}(\check{t},\check{s}) = \false$. 
According to \cref{shifted_support_bbox} and \cref{A_eta}, we therefore get 
\begin{align*}
& \max\{\diam(B^{(\ell-1)}), \diam(C^{(\ell-1)})\}
= \max\{\diam(B_{\check{t}}), \diam(C_{\check{s}})\}
> \eta\dist(B_{\check{t}},C_{\check{s}})\\
&= \eta\inf\{\|x + m_{\check{t}} - (y + m_{\check{s}})\|_{2} : x \in B^{(\ell-1)}, y \in C^{(\ell-1)}\}\\
&= \eta\inf\big\{\big\|x - \big(y - (m_{\check{t}} - m_{\check{s}})\big)\big\|_{2} : x \in B^{(\ell-1)}, y \in C^{(\ell-1)}\big\}\\
&= \eta\dist\big(B^{(\ell-1)},C^{(\ell-1)} - (m_{\check{t}} - m_{\check{s}})\big),
\end{align*}
which, in view of \cref{all_placements}, \cref{general_placement},
\cref{shifted_ref_supp_bbox} and \cref{Ahat_eta}, implies 
\begin{equation}\label{ref_inadm}
r^{(\ell-1)} - (m_{\check{t}} - m_{\check{s}}) = r^{(\ell-1)} + \delta^{(\ell-1)} \odot (p_{\check{s}} - p_{\check{t}}) \in \underbrace{\big\{u \in \that^{(\ell-1)} : \ahat_{\eta}(r^{(\ell-1)},u) = \false\big\}}_{=:H^{(\ell-1)}}.
\end{equation}
Moreover, \cref{sondiffvec} implies
\begin{align*}
m_{t} - m_{s} \in 
\{m_{\check{t}} - m_{\check{s}},
 m_{\check{t}} - m_{\check{s}} + \delta^{(\ell)}_{\iota^{(\ell-1)}}e_{\iota^{(\ell-1)}},
 m_{\check{t}} - m_{\check{s}} - \delta^{(\ell)}_{\iota^{(\ell-1)}}e_{\iota^{(\ell-1)}}\},
\end{align*}
which, combined with \cref{ref_inadm} and \cref{general_placement}, leads to
\begin{equation*}
m_{t} - m_{s} \in 
\bigcup_{u \in H^{(\ell-1)}}\{-m_{u}, -m_{u} + \delta^{(\ell)}_{\iota^{(\ell-1)}}e_{\iota^{(\ell-1)}}, 
-m_{u} - \delta^{(\ell)}_{\iota^{(\ell-1)}}e_{\iota^{(\ell-1)}}\}.
\end{equation*}
Since we have chosen $ (t,s) \in \ctIJ^{(\ell)} $ arbitrarily, this leads to
\begin{equation*}
|\{m_{t} - m_{s} : (t,s) \in \ctIJ^{(\ell)}\}| \leq 3|H^{(\ell-1)}|.
\end{equation*}
Furthermore, due to \cref{general sparsity}, we have
$ |H^{(\ell-1)}| \leq \const_{sp} $. \hfill $ \Box $

%
%
\begin{theorem}[Storage requirements of the coupling matrices]
\label{coupling storage}
The coupling matrices $ (S_{b})_{b\in\Leaves_{\ctIJ}^{+}} $ require not
more than $3\const_{sp}k^{2}(\ell_{max}-d)$ units of storage.
\end{theorem}
\textit{Proof.} On every level $\ell \in \{0,...,d\}$ 
all boxes $ t \in \ctI^{(\ell)} $ and $ s \in \ctJ^{(\ell)} $ contain, by
construction, the midpoint of $ r^{(0)} $, which implies
$\dist(B_{t},C_{s}) = 0$ for all $ (t,s) \in \ctIJ^{(\ell)}$
and therefore
$\Leaves_{\ctIJ}^{+}\cap\ctIJ^{(\ell)}
= \{(t,s) \in \ctIJ^{(\ell)} : A_{\eta}(t,s) = \true\} = \emptyset$.
Combining this with the implication \cref{tinv_coupling}, \cref{diff_bound}
and \cref{nk} directly implies that for representing all coupling matrices
$(S_{b})_{b\in\Leaves_{\ctIJ}^{+}}$ not more than
\begin{equation*}
|\Theta|^{2}\sum_{\ell = d+1}^{\ell_{max}}\left|\{m_{t} - m_{s} : (t,s) \in \Leaves_{\ctIJ}^{+}\cap\ctIJ^{(\ell)}\}\right|
\leq
k^{2}\sum_{\ell = d+1}^{\ell_{max}}3\const_{sp}
=
3\const_{sp}k^{2}(\ell_{max}-d)
\end{equation*}
units of storage are needed. $\hfill \Box$

In order for arithmetic operations like the matrix-vector multiplication
to have the desired complexity, $ \ctIJ $ has to be \textit{sparse}.
This property is guaranteed by our construction.

%
%
\begin{lemma}[Sparsity of $ \ctIJ $]\label{blocktree_sparsity}
We have
\begin{align*}
\left|\{s \in \ctJ : (t,s) \in \ctIJ\}\right| &\leq 2\const_{sp}
\qquad \text{for all $ t \in \ctI $}, \\
\left|\{t \in \ctI : (t,s) \in \ctIJ\}\right| &\leq 2\const_{sp}
\qquad \text{for all $ s \in \ctJ $}.
\end{align*}
\end{lemma}
\textit{Proof.} 
The statements are implied by the definition of $ \ctIJ $,
\cref{all same level}, \cref{general sparsity} and the fact that
every node of $\ctI$ or $\ctJ$ has at most two sons. $\hfill\Box$

For the sake of completeness, we also have to consider the storage
requirements of the \emph{nearfield matrices}
\begin{align*}
  N_{(t,s)} &:= G|_{I_t\times J_s} &
  &\text{ for all } (t,s) \in \Leaves_{\ctIJ}^{-}
\end{align*}
that capture the part of the matrix that cannot be compressed.
Since we want to maintain the ability to work with unstructured
surface meshes, we cannot avoid storing all of these matrices explicitly.

%
%
\begin{lemma}[Storage requirements of the nearfield matrices]
\label{nearfield storage}
We define the \emph{resolutions} of the cluster trees $\ctI$ and $\ctJ$ by
\begin{align*}
  \gamma_{I} &:= \max\{|I_{t}| : t \in \Leaves_{\ctI}\}, &
  \gamma_{J} &:= \max\{|J_{s}| : s \in \Leaves_{\ctJ}\}.             
\end{align*}
Storing $(N_{b})_{b \in \Leaves_{\ctIJ}^{-}}$  requires not more than
$\const_{sp}\min\{\gamma_{I},\gamma_{J}\}n $ units of storage.
\end{lemma}
\noindent\textit{Proof.} 
By construction and due to (\ref{leaves_lmax}), (\ref{translational}) and (\ref{all_placements}) we have
\begin{equation*}
\Leaves_{\ctIJ}^{-} \, \subseteq \, \Leaves_{\ctI} \times \Leaves_{\ctJ} \, = \, \ctI^{(\ell_{max})} \times \ctJ^{(\ell_{max})} \, \subseteq \, \that^{(\ell_{max})} \times \that^{(\ell_{max})}.
\end{equation*}
Lemma \ref{general sparsity} and (\ref{leafpart}) therefore imply that 
storing all nearfield matrices $ (\Ngen_{b})_{b\in\Leaves_{\ctIJ}^{-}} $ requires
not more than
\begin{align*}
\sum_{(t,s) \in \Leaves_{\ctIJ}^{-}}|I_{t}|\,|J_{s}|
&= \sum_{t\in \Leaves_{\ctI}}\sum_{\substack{s\in \Leaves_{\ctJ}\\(t,s)\in\Leaves_{\ctIJ}^{-}}}|I_{t}|\,|J_{s}|
\leq \sum_{t\in \Leaves_{\ctI}}\sum_{\substack{s\in \that^{(\ell_{max})}\\\ahat_{\eta}(t,s) = \false}}\gamma_{J}|I_{t}|\\
&\leq \sum_{t\in \Leaves_{\ctI}}\const_{sp}\gamma_{J}|I_{t}|
= \const_{sp}\gamma_{J}\sum_{t\in \Leaves_{\ctI}}|I_{t}|
= \const_{sp}\gamma_{J}|I| 
\end{align*}
units of storage and similarly not more than $ \const_{sp}\gamma_{I}|J| $ units of storage, which, 
in view of (\ref{nk}),
implies the statement. $\hfill\Box$

\Cref{transfer storage}, \cref{nearfield storage} and \cref{coupling storage}
suggest that we still have to discuss how to properly choose the maximal
level $\ell_{max}$.
On the one hand, we have to choose $ \ell_{max} $ large enough to keep at
least one of the resolutions
$ \gamma_{I} $ and $ \gamma_{J} $
sufficiently small.
On the other hand, since the complexity of the matrix-vector
multiplication mainly depends on $|\ctI|$ and $|\ctJ| $
(cf. \cite[Theorem 3.42]{BO10}), we have to choose $\ell_{max}$ small enough
to prevent $|\ctI|$ and $|\ctJ|$ from getting too large.

In order to be able to analyze $ \gamma_{I} $ and $ \gamma_{J} $, we assume
that the relative number of characteristic points assigned to a box
(cf. (\ref{index_subsets}) and (\ref{index_subsets_2})) can essentially be estimated by its size,
i.e., we assume that there is a constant $ \const_{nd} \in \bbbr_{>0}$
satisfying
\begin{equation}
\label{number_diameter}
|I_{t}| \leq \const_{nd}\diam(t)^{d-1}|I|,
\quad
|J_{s}| \leq \const_{nd}\diam(s)^{d-1}|J|
\qquad\text{for all $ t \in \ctI, s \in \ctJ $}.
\end{equation} 
This condition is, for example, usually fulfilled on \emph{shape-regular} and \emph{quasi-uniform} meshes (cf. \cite[Remark 4.1.14]{SASC11}).
Due to (\ref{regular_splitting}) and (\ref{halved_ref_box}), we have
\begin{equation}\label{halved_diameter}
\diam(r^{(\ell+d)}) = \tfrac{1}{2}\diam(r^{(\ell)}) \qquad \text{for all $ \ell \in \{0,...,\ell_{max}-d\}$}.
\end{equation}
For the sake of simplicity, we therefore choose $ \ell_{max} \in d\cdot\bbbn $.
\noindent According to (\ref{number_diameter}), (\ref{translational}) and
(\ref{halved_diameter}), we thus get
\begin{align*}
|I_{t}| &\leq \const_{nd} \diam(t)^{d-1}|I| \\
&= \const_{nd}\diam(r^{(\ell_{max})})^{d-1}|I|\\
&= \const_{nd}2^{-\frac{d-1}{d}\ell_{max}}\diam(r^{(0)})^{d-1}|I|
  \qquad\text{ for all } t\in\ctI^{(\ell_{max})} = \Leaves_{\ctI}
\end{align*}
and likewise
\begin{align*}\label{res_TJ}
|J_{s}| 
&\leq 
\const_{nd}2^{-\frac{d-1}{d}\ell_{max}}\diam(r^{(0)})^{d-1}|J|
  \qquad\text{ for all } s \in \ctJ^{(\ell_{max})} = \Leaves_{\ctJ}.
\end{align*}
Hence, with
$ \const_{rs} 
:= \const_{nd}\diam(r^{(0)})^{d-1} $
we have
\begin{equation}\label{resolutions}
\gamma_{I} \leq \const_{rs}2^{-\frac{d-1}{d}\ell_{max}}|I|,
\qquad
\gamma_{J} \leq \const_{rs}2^{-\frac{d-1}{d}\ell_{max}}|J|.
\end{equation}
We choose a constant $ \const_{rk} \in \bbbr_{>0} $ satisfying $\const_{rk}k \leq \min\{|I|, |J|\}$ and observe 
\begin{align*}
2^{-\frac{d-1}{d}\ell_{max}}|I| \leq \const_{rk}k 
&\iff \ell_{max} \geq \tfrac{d}{d-1}\big(\log_{2}(|I|) - \log_{2}(\const_{rk}k)\big),\\
2^{-\frac{d-1}{d}\ell_{max}}|J| \leq \const_{rk}k &\iff \ell_{max} \geq \tfrac{d}{d-1}\big(\log_{2}(|J|) - \log_{2}(\const_{rk}k)\big).
\end{align*}
Based on this, we choose 
\begin{equation}\label{ellmax}
 \ell_{max} := \min\big\{\ell \in d \cdot \bbbn : \ell \geq
 \tfrac{d}{d-1}\big(\min\{\log_{2}(|I|), \log_{2}(|J|)\} - \log_{2}(\const_{rk}k)\big)\big\}
\end{equation}
to gain, according to \cref{resolutions},
\begin{equation}\label{res_bound}
\min\{\gamma_{I}, \gamma_{J}\} \leq \const_{rs}\const_{rk}k.
\end{equation}

In order to be able to analyze $ |\ctI| $ and $ |\ctJ| $, we assume that the
number of boxes per level which contain only a relatively small part of the
boundary $ \Gamma $ does not grow too fast, i.e., we assume that there are
constants $ \const_{da}, \const_{sa} \in \bbbr_{>0} $ satisfying
\begin{subequations}\label{diameter_area}
\begin{align}
|\{t \in \ctI^{(\ell)}: \diam(t)^{d-1} > \const_{da}|t\cap\Gamma|\}| &\leq \const_{sa}\diam(r^{(\ell)})^{-(d-1)}\\
|\{s \in \ctJ^{(\ell)}: \diam(s)^{d-1} > \const_{da}|s\cap\Gamma|\}| &\leq \const_{sa}\diam(r^{(\ell)})^{-(d-1)}
\end{align}
\end{subequations}
on every level $ \ell \in \{0,...,\ell_{max}\} $, where 
\begin{equation*}
|X| := \int_{\Gamma}\mathbbm{1}_{X}(y) \, dy \qquad \text{for all measurable subsets $ X \subseteq \Gamma$}.
\end{equation*}
\begin{lemma}\label{number_of_big_area_boxes}
	On every level $ \ell \in \{0,...,\ell_{max}\} $ we have
	\begin{align*}
	|\{t \in \ctI^{(\ell)} : \diam(t)^{d-1} \leq \const_{da}|t\cap\Gamma|\}| \,\leq\, 2\const_{da}|\Gamma|\diam(r^{(\ell)})^{-(d-1)},\\
	|\{s \in \ctJ^{(\ell)} : \diam(s)^{d-1} \leq \const_{da}|s\cap\Gamma|\}| \,\leq\, 2\const_{da}|\Gamma|\diam(r^{(\ell)})^{-(d-1)}.
	\end{align*}
\end{lemma}
\textit{Proof.}
\noindent Let $ \ell \in \{0,...,\ell_{max}\} $ and let
$ T_{I,\const_{da}}^{(\ell)} := \{t \in \ctI^{(\ell)} : \diam(t)^{d-1} \leq \const_{da}|t\cap\Gamma|\} $. 
Due to (\ref{translational}), the surface area of the set
$ \big\{x \in \Gamma : |\{t \in T^{(\ell)}_{I,\const_{da}} : x \in t\}| > 2\big\} $ equals $ 0 $. Thus
\begin{align*}
\diam(r^{(\ell)})^{d-1}\Big|T^{(\ell)}_{I,\const_{da}}\Big| &= \sum_{t \in T^{(\ell)}_{I,\const_{da}}}\diam(t)^{d-1}
\leq \sum_{t\in T^{(\ell)}_{I,\const_{da}}}\const_{da}|t\cap\Gamma|\\
&= \const_{da}\sum_{t\in T^{(\ell)}_{I,\const_{da}}}\int_{\Gamma}\mathbbm{1}_{t\cap\Gamma}(y)\,dy
= \const_{da}\int_{\Gamma}\sum_{t\in T^{(\ell)}_{I,\const_{da}}}\mathbbm{1}_{t\cap\Gamma}(y)\,dy\\
&\leq \const_{da}\int_{\Gamma}2\cdot\mathbbm{1}_{\Gamma}(y)\,dy
= 2\const_{da}|\Gamma|
\end{align*}
holds, which proves the first statement. The second statement can be proven similarly. $ \hfill \Box $\\

Defining $ \const_{nb} := 2\const_{da}|\Gamma| + \const_{sa} $ and combining
assumption (\ref{diameter_area}) with Lemma \ref{number_of_big_area_boxes}
yield
\begin{subequations}\label{number_of_boxes}
\begin{align}
|\ctI^{(\ell)}| &\leq \const_{nb}\diam(r^{(\ell)})^{-(d-1)} &
  &\text{ for all } \ell \in \{0,...,\ell_{max}\},\\
|\ctJ^{(\ell)}| &\leq \const_{nb}\diam(r^{(\ell)})^{-(d-1)} &
  &\text{ for all } \ell \in \{0,...,\ell_{max}\}.
\end{align} 
\end{subequations}
Moreover, according to (\ref{ellmax}), we also have 
\begin{equation}\label{lmax_bound}
\ell_{max} \leq 
\tfrac{d}{d-1}\big(\min\{\log_{2}(|I|), \log_{2}(|J|)\} - \log_{2}(\const_{rk}k)\big) + d.
\end{equation}
By using (\ref{number_of_boxes}), (\ref{halved_diameter}) and (\ref{lmax_bound})
we get
\begin{align*}
|\ctI| &= \sum_{\ell=0}^{\ell_{max}}|\ctI^{(\ell)}|
\leq \sum_{\ell=0}^{\ell_{max}}\const_{nb}\diam(r^{(\ell)})^{-(d-1)}
= \const_{nb}\sum_{\ell=0}^{\ell_{max}}\diam(r^{(\ell)})^{-(d-1)}\\
&= \const_{nb}\Bigg(\diam(r^{(0)})^{-(d-1)} + 
\sum_{\ell=0}^{\frac{\ell_{max}}{d}-1}\sum_{u = 1}^{d}\diam(r^{(d\ell+u)})^{-(d-1)}\Bigg)\\
&<
\const_{nb}\Bigg(\diam(r^{(0)})^{-(d-1)} + 
\sum_{\ell=0}^{\frac{\ell_{max}}{d}-1}d\diam(r^{(d(\ell+1))})^{-(d-1)} \Bigg)\\
&=
\const_{nb}\Bigg(\diam(r^{(0)})^{-(d-1)} + 
d\sum_{\ell=0}^{\frac{\ell_{max}}{d}-1}2^{(\ell+1)(d-1)}\diam(r^{(0)})^{-(d-1)} \Bigg)\\
&<
\const_{nb}\Bigg(d\diam(r^{(0)})^{-(d-1)} + 
d\diam(r^{(0)})^{-(d-1)}\sum_{\ell=0}^{\frac{\ell_{max}}{d}-1}2^{(\ell+1)(d-1)} \Bigg)\\
&=
\const_{nb}d\diam(r^{(0)})^{-(d-1)}\Bigg(1 + 
\sum_{\ell=1}^{\frac{\ell_{max}}{d}}\big(2^{d-1}\big)^{\ell} \Bigg)\\
&=
\const_{nb}d\diam(r^{(0)})^{-(d-1)}
\sum_{\ell=0}^{\frac{\ell_{max}}{d}}\big(2^{d-1}\big)^{\ell}
=
\const_{nb}d\diam(r^{(0)})^{-(d-1)}\frac{1-(2^{d-1})^{\frac{\ell_{max}}{d}+1}}{1 - 2^{d-1}}\\
&<
\const_{nb}d\diam(r^{(0)})^{-(d-1)}(2^{d-1}-1)^{-1}2^{\frac{d-1}{d}\ell_{max} + d-1}\\
&\leq
\const_{nb}d\diam(r^{(0)})^{-(d-1)}(2^{d-1}-1)^{-1}2^{\frac{d-1}{d}(\frac{d}{d-1}(\log_{2}(|I|) - \log_{2}(\const_{rk}k)) + d) + d-1}\\
&=
\const_{nb}d\diam(r^{(0)})^{-(d-1)}(2^{d-1}-1)^{-1}2^{\log_{2}(|I|) - \log_{2}(\const_{rk}k) + 2(d-1)}\\
&=
\const_{nb}d\diam(r^{(0)})^{-(d-1)}(2^{d-1}-1)^{-1}2^{2(d-1)}\const_{rk}^{-1}k^{-1}|I|.
\end{align*}
Due to $ \frac{x^{2}}{x-1} = \frac{x^{2}(x+2)}{(x-1)(x+2)} = \frac{x^{2}(x+2)}{x^{2}+x-2}\leq x+2 $ for $ x \in \bbbr_{\geq2} $ we therefore have
\begin{equation*}
|\ctI| < \const_{nb}d\diam(r^{(0)})^{-(d-1)}(2^{d-1} + 2)\const_{rk}^{-1}k^{-1}|I|
\end{equation*}
and analogously
\begin{equation*}
|\ctJ| < \const_{nb}d\diam(r^{(0)})^{-(d-1)}(2^{d-1}+2)\const_{rk}^{-1}k^{-1}|J|.
\end{equation*}

Finally, based on our choice of the maximal level $ \ell_{max} $, we can now specify the storage requirements for our $ \htwo$-matrix approximation:
\begin{theorem}[Total storage requirements]\label{total_storage}
	$ \widetilde{G} $ requires not more than
	\begin{itemize}
		\item $ 2kn $ units of storage for the leaf matrices,
		\item 
		$ 4k^{2}\big(\frac{d}{d-1}\big(\log_{2}(n) - \log_{2}(\const_{rk}k)\big) + d\big)  $
		units of storage 
        for the transfer matrices,
        \item $ \const_{sp}\const_{rs}\const_{rk}kn $ units of storage for the nearfield matrices,
        \item $ 3\const_{sp}k^{2}\frac{d}{d-1}\big(\log_{2}(n) - \log_{2}(\const_{rk}k)\big) $
        units of storage 
        for the coupling matrices.
	\end{itemize}
\end{theorem}
\textit{Proof.} The statement is a direct consequence of Lemma \ref{leaf storage}, Theorem \ref{transfer storage}, Lemma \ref{nearfield storage}, Theorem \ref{coupling storage}, (\ref{res_bound}) and (\ref{lmax_bound}). \hfill $ \Box $

Compared to the standard $ \htwo $-approach, these storage complexity bounds
lead to a significant reduction of the total storage requirements, as we will
see and discuss in chapter 4.

%
%
\begin{remark}[Bounded enlargement]\label{assumption_true}
Due to (\ref{halved_diameter}) and (\ref{lmax_bound}), we have
\begin{equation*}
\diam(r^{(\ell)}) \geq \diam({r^{(\ell_{max})}}) \sim 2^{-\frac{\ell_{max}}{d}} \gtrsim 2^{-\frac{\log_{2}(n)}{d-1}} =  n^{-\frac{1}{d-1}} \qquad \text{for all $ \ell \in \{0,...,\ell_{max}\} $}.
\end{equation*}	
In many situations, the diameters of the supports behave like
$n^{-\frac{1}{d-1}}$, for example on shape-regular and quasi-uniform meshes.
In these cases, the condition (\ref{diameter_outer_inner}) is therefore fulfilled.
\end{remark}

\section{Numerical experiments}
\label{se:experiments}

In the following, let 
\begin{equation*}
g: \bbbr^{d} \times \bbbr^{d} \rightarrow \bbbr, \quad
(x,y) \mapsto  
\begin{cases}
\frac{1}{4\pi}\frac{1}{\|x-y\|_{2}} &\text{if $ x \neq y$}\\
0                                   &\text{if $ x = y $}
\end{cases}
\end{equation*}
and for all $ y \in \Gamma = \partial \Omega$ let $ n_{y} $ denote the corresponding outward-pointing unit normal vector. Under suitable conditions, a solution $ u $ of the \emph{interior Dirichlet boundary value problem}
\begin{equation}\label{laplace_hom_dir}
\begin{split}
-\Delta u &= 0 \qquad\text{in $\Omega $},\\
u &= g_{D} \quad\text{ on $\Gamma $}
\end{split}
\end{equation}
is given by \emph{Green's representation formula} (cf. \cite[Theorem 2.2.2]{HA92}, \cite[Theorem 3.1.6]{SASC11}), which reduces obtaining $ u $ to obtaining the corresponding \emph{Neumann boundary values} $ g_{N} : \Gamma \rightarrow \bbbr, y \mapsto \frac{\partial u}{\partial n_{y}}(y)$. Depending on the situation, $ g_{N} $ as well as the given \emph{Dirichlet boundary values} $ g_{D} $ might only exist in a generalized sense. The Neumann values $ g_{N} $ can be obtained by solving a boundary integral equation resulting from combining Green's representation formula with suitable trace operators (cf. \cite[Section 3.4.2.1]{SASC11}).

Applying a Galerkin discretization to this integral equation on a triangulation of $ \Gamma $ with discontinuous piecewise constant basis functions $ (\varphi_{i})_{i\in I} $ and replacing $ g_{D} $ with the corresponding $ L^{2} $-orthogonal projection $ \tilde{g}_{D} $ into a space spanned by continuous piecewise linear nodal basis functions $ (\psi_{j})_{j\in J} $ leads to a linear system
\begin{equation}\label{Galerkin_LS}
Gx = \Big(K + \frac{1}{2}M\Big)b,
\end{equation}
where $ b \in \bbbr^{J} $ contains the coefficients of $ \tilde{g}_{D} = \sum_{j \in J}b_{j}\psi_{j} $, $ x \in \bbbr^{I} $ contains the coefficients of the approximate Neumann values $ \tilde{g}_{N} := \sum_{i \in I}x_{i}\varphi_{i} $, and where $ G \in \bbbr^{I\times I} $, $ K \in \bbbr^{I \times J} $ and $ M \in \bbbr^{I \times J} $ are given by
\begin{align*}
G_{ij} &:= \int_\Gamma \varphi_i(x)
\int_\Gamma g(x,y) \varphi_j(y) \,dy\,dx 
&&\text{ for all } i \in I,\ j \in I,\\
K_{ij} &:= \int_\Gamma \varphi_i(x)
\int_\Gamma \frac{\partial g}{\partial n_{y}}(x,y)\, \psi_{j}(y) \,dy\,dx 
&&\text{ for all } i \in I,\ j \in J,\\
M_{ij} &:= \int_\Gamma \varphi_i(x)
\psi_j(x) \,dx
&&\text{ for all } i \in I,\ j \in J.
\end{align*}
The derivative $ \frac{\partial}{\partial n_{y}} $ is applied with respect to the $ y $-variable.

We construct an $ \htwo $-matrix approximation $ \widetilde{G} $ of $ G $ as described in \cref{se:h2matrices}. In order to approximate $ K $, we replace $ \frac{\partial}{\partial n_{y}}g(x,y) $ with $ \frac{\partial}{\partial n_{y}}(\interpol_{B_{t}} \otimes \interpol_{C_{s}})[g](x,y) $ for all $(x,y) \in B_{t} \times C_{s}$, $ (t,s) \in \Leaves_{\ctIJ}^{+} $ (cf. \cref{tensor_approx}). This approach leads to an $ \htwo $-matrix approximation $ \widetilde{K} $ of $ K $ which differs from the $ \htwo $-matrix described in \cref{se:h2matrices} only in terms of the leaf matrices for the column cluster basis, which are given by
\begin{equation*}
W_{s} \in \bbbr^{J_{s} \times \Theta}, \quad
(W_{s})_{j\mu} := \int_{\Gamma}\psi_{j}(y)\frac{\partial}{\partial n_{y}}\lagrange_{C_{s}, \mu}(y) \, dy
\qquad \text{for all $ s \in \Leaves_{\ctJ}, j \in J_{s}, \mu \in \Theta $}.
\end{equation*}
Hence, all results from \cref{se:complexity} also hold for $ \widetilde{K} $. Moreover, we can expect a similar convergence behavior as for $ \widetilde{G} $ (cf. \cite[Chapter 4]{BO10}). The sparse matrix $ \frac{1}{2}M $ can simply be added to the nearfield entries of $ \widetilde{K} $ and therefore does not require any additional storage. 

In this way, we approximate the Neumann values of the function $ u_{0} := 4\pi g( \cdot , y_{0}) $ with $ y_{0} := (1.2, 1.2, 1.2)^{T} $ on an aproximation $ \Gamma_{S} $ of the unit sphere \small$S^2 := \{x \in \bbbr^{3}:\|x\|_{2} = 1\}$ 
= $\partial\{x \in \bbbr^{3}:\|x\|_{2} < 1\}$ 
\normalsize con\-sisting of plane triangles. $ S^2 $ and its approximation are both contained in the box $ [-1,1]^{3} $, which we use as our starting box $ r^{(0)} $. We choose the centers of gravity of the triangles as characteristic points for $ (\varphi_{i})_{i\in I} $, the vertices of the triangles as characteristic points for $ (\psi_{j})_{j\in J}, \const_{rk} = 2$ and $ \eta = 2 $. We denote the maximal level of the trees used for the construction of $ \widetilde{G} $ by $ \ell_{max}^{I \times I} $ and the maximal level of the trees used for the construction of $ \widetilde{K} $ by $ \ell_{max}^{I\times J} $. We compute the matrix entries by tensor Gauss quadrature as described in \cite{ERSA98, SA96},\cite[Chapter 5]{SASC11} with $ 3 $ quadrature points per dimension for the regular integrals and $ 5 $ quadrature points per dimension for the singular integrals. In order to preserve the convergence behavior of the standard Galerkin solution, we increase the interpolation degree $ \theta $ by $ 1 $ whenever the mesh width $ h $ is halved, i.e., whenever the number $ n = \max\{|I|,|J|\} = |I|$ of triangles is quadrupled (cf. \cite[Section 10.1 and Section 10.2]{BO10}). We solve the linear system using the cg method with a relative residual accuracy of $10^{-6}$ (cf. \cite[Section 6.1]{SASC11}).

\begin{table}[H]
	\centering
	\begin{tabular}{ |r|r|r|r|r|r| }
		$ n $  & $ \theta $ & $ k$ &  $ \ell_{max}^{I\times I} $ & $ \ell_{max}^{I\times J} $ &
		$ \epsilon_{L^{2}} $ \\ \hline 
		8192 & 4 & 125& 9 & 9 & $ 1.194_{-2}$\\
		32768 & 5 & 216& 12 & 9 & $ 5.683_{-3}$\\
		131072 & 6 & 343 & 12 & 12 & $ 2.870_{-3}$\\
		524288 & 7 & 512 & 15 & 15 & $ 1.412_{-3}$\\
		2097152 & 8 & 729 & 18 & 15 & $ 7.050_{-4}$ 
	\end{tabular}
	\caption{Parameters and resulting $ L^{2} $-error $\epsilon_{L^{2}} := \|g_{N} - \tilde{g}_{N}\|_{L^{2}(\Gamma_{S})}$.}
	\label{par_err_tab}
\end{table}

\begin{table}[H]
	\centering
	\begin{tabular}{ |r|r|r|r|r| }
		$ n $  & \small leaf matrices & \small transfer matrices & \small nearfield matrices & \small coupling matrices\\\hline  
		8192 & 15.6 & 4.3 & 74.4 &51.8 \\
		32768 & 108.1 & 17.1 & 303.8 &218.6 \\
		131072 & 686.1 & 43.1 & 3913.3 &509.9 \\
		524288 & 4096.3  & 120.0 & 16090.5 &1496.1 \\
		2097152 & 23329.2 & 291.9 & 65051.2 &3762.8  
	\end{tabular}
	\caption{Storage requirements of $ \widetilde{G} $ on $ \Gamma_{S} $ in MB.}
	\label{storage_Gtilde_tab}
\end{table}

\begin{table}[H]
	\centering
	\begin{tabular}{ |r|r|r|r|r| }
		$ n $  & \small leaf matrices & \small transfer matrices & \small nearfield matrices & \small coupling matrices \\ \hline
		8192 & 11.7 & 4.3 & 38.8 &54.2  \\
		32768 & 81.0 & 12.8 & 476.1 &143.9 \\
		131072 & 514.6 & 43.1 & 1956.3 &524.3  \\
		524288 & 3072.3 & 120.0 & 8047.3 &1528.1  \\
		2097152 & 17496.3 & 243.3  & 128696.6 &2765.3   
	\end{tabular}
	\caption{Storage requirements of $ \widetilde{K} $ on $ \Gamma_{S} $ in MB.}
	\label{storage_Ktilde_tab}
\end{table}

It is clearly visible in the Tables \ref{par_err_tab}, \ref{storage_Gtilde_tab} and \ref{storage_Ktilde_tab} that the storage requirements of our $ \htwo $-matrices $ \widetilde{G} $ and $ \widetilde{K} $ behave as predicted in Chapter 3: The storage for the leaf matrices grows like $ \mathcal{O}(nk) $ (cf. Lemma \ref{leaf storage}), the storage for the transfer matrices grows like $ \mathcal{O}(\ell_{max}k^{2}) $ (cf. Theorem \ref{transfer storage}), the storage for the nearfield matrices grows like $ \mathcal{O}(2^{-\frac{2}{3}\ell_{max}}n^{2}) $ (cf. Lemma \ref{nearfield storage}, (\ref{resolutions})) and the storage for the coupling matrices grows like $  \mathcal{O}((\ell_{max}-3)k^{2}) $ (cf. Theorem \ref{coupling storage}).

%
%
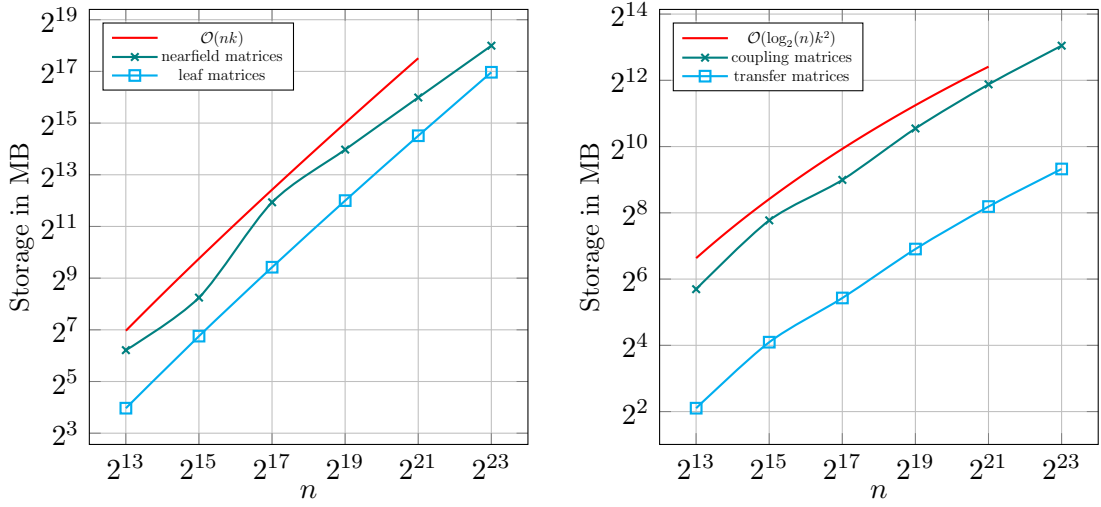
\begin{figure}[H]
\begin{tikzpicture}
\begin{axis}[
xmode = log,
log basis x={2},
ymode = log,
log basis y={2},
xtick distance = 4.0,  
ytick distance = 4.0,
grid = both,
width = 0.5\textwidth,
height = 0.5\textwidth,
xlabel = {$ n $},
ylabel = {\small Storage in MB},
x label style={at={(axis description cs:0.5,0.02)}},
y label style={at={(axis description cs:0.06,.5)}},
legend pos = north west,
legend style={nodes={scale=0.5}}]
\newcommand\m{(log2(x)/2.0) - 2.5};
\newcommand\rank{(\m + 1)^3};
\addplot[smooth, thick, red, domain = 8192:2097152]{16 * x *\rank * 8.0/(1024.0 * 1024.0)};
\addplot[smooth, thick, teal, mark = x] table[x={n}, y={nearfield_storage}] {storage_leaf_nearfield_Gtilde.dtn};
\addplot[smooth, thick, cyan, mark = cube] table[x={n}, y={leaf_storage}] {storage_leaf_nearfield_Gtilde.dtn};
\legend{\large$ \mathcal{O}(nk) $, \large nearfield matrices, \large leaf matrices}
\end{axis}

\begin{axis}[xshift = 7.5cm,
xmode = log,
log basis x={2},
ymode = log,
log basis y={2},
xtick distance = 4.0,  
ytick distance = 4.0,
grid = both,
width = 0.5\textwidth,
height = 0.5\textwidth,
xlabel = {$ n $},
ylabel = {\small Storage in MB},
x label style={at={(axis description cs:0.5,0.02)}},
y label style={at={(axis description cs:0.06,.5)}},
legend pos = north west,
legend style={nodes={scale=0.5}}]
\newcommand\m{(log2(x)/2.0) - 2.5};
\newcommand\rank{(\m + 1)^3};
\addplot[smooth, thick, red, domain = 8192:2097152]{64 * log2(x) *\rank^2 * 8.0/(1024.0 * 1024.0)};
\addplot[smooth, thick, teal, mark = x] table[x={n}, y={coupling_storage}] {storage_transfer_coupling_Gtilde.dtn};
\addplot[smooth, thick, cyan, mark = cube] table[x={n}, y={transfer_storage}] {storage_transfer_coupling_Gtilde.dtn};
\legend{\large$ \mathcal{O}(\log_{2}(n)k^{2}) $, \large coupling matrices, \large transfer matrices}
\end{axis}
\end{tikzpicture}
\caption{Storage requirements of $ \widetilde{G} $ on $ \Gamma_{S} $ in MB.}
\label{storage_Gtilde_fig}
\end{figure}

\begin{figure}[H]
	\begin{tikzpicture}
	\begin{axis}[
	xmode = log,
	log basis x={2},
	ymode = log,
	log basis y={2},
	xtick distance = 4.0,  
	ytick distance = 4.0,
	grid = both,
	width = 0.5\textwidth,
	height = 0.5\textwidth,
	xlabel = {$ n $},
	ylabel = {\small Storage in MB},
	x label style={at={(axis description cs:0.5,0.02)}},
	y label style={at={(axis description cs:0.06,.5)}},
	legend pos = north west,
	legend style={nodes={scale=0.5}}]
	\newcommand\m{(log2(x)/2.0) - 2.5};
	\newcommand\rank{(\m + 1)^3};
	\addplot[smooth, thick, red, domain = 8192:2097152]{16 * x *\rank * 8.0/(1024.0 * 1024.0)};
	\addplot[smooth, thick, teal, mark = x] table[x={n}, y={nearfield_storage}] {storage_leaf_nearfield_Ktilde.dtn};
	\addplot[smooth, thick, cyan, mark = cube] table[x={n}, y={leaf_storage}] {storage_leaf_nearfield_Ktilde.dtn};
	\legend{\large$ \mathcal{O}(nk) $, \large nearfield matrices, \large leaf matrices}
	\end{axis}
	
	\begin{axis}[xshift = 7.5cm,
	xmode = log,
	log basis x={2},
	ymode = log,
	log basis y={2},
	xtick distance = 4.0,  
	ytick distance = 4.0,
	grid = both,
	width = 0.5\textwidth,
	height = 0.5\textwidth,
	xlabel = {$ n $},
	ylabel = {\small Storage in MB},
	x label style={at={(axis description cs:0.5,0.02)}},
	y label style={at={(axis description cs:0.06,.5)}},
	legend pos = north west,
	legend style={nodes={scale=0.5}}]
	\newcommand\m{(log2(x)/2.0) - 2.5};
	\newcommand\rank{(\m + 1)^3};
	\addplot[smooth, thick, red, domain = 8192:2097152]{64 * log2(x) *\rank^2 * 8.0/(1024.0 * 1024.0)};
	\addplot[smooth, thick, teal, mark = x] table[x={n}, y={coupling_storage}] {storage_transfer_coupling_Ktilde.dtn};
	\addplot[smooth, thick, cyan, mark = cube] table[x={n}, y={transfer_storage}] {storage_transfer_coupling_Ktilde.dtn};
	\legend{\large$ \mathcal{O}(\log_{2}(n)k^{2}) $, \large coupling matrices, \large transfer matrices}
	\end{axis}
	\end{tikzpicture}
	\caption{Storage requirements of $ \widetilde{K} $ on $ \Gamma_{S} $ in MB.}
	\label{storage_Ktilde_fig}
\end{figure}
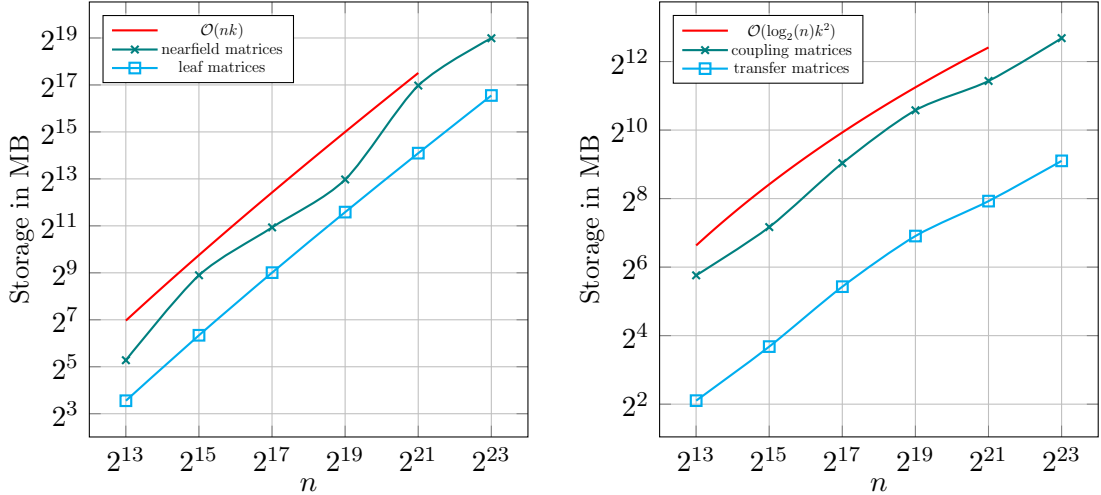

Figure \ref{storage_Gtilde_fig} and Figure \ref{storage_Ktilde_fig} show that the storage requirements of the matrix components can be bounded as predicted by Theorem \ref{total_storage}. The expected theoretical convergence rate of $ \mathcal{O}(h) $ for the Neumann values in the $ L^{2}$-norm (cf. \cite[Theorem 4.1.33]{SASC11} and \cite[Theorem 4.6]{DAFAGRHASA01}) is preserved (see Table \ref{par_err_tab}). 

We want to compare these results with the corresponding results obtained by
conventional $ \htwo $-matrix approximations $ \widetilde{G}_{conv} $ and $
\widetilde{K}_{conv} $, which are based on different cluster trees and support
bounding boxes (cf. \cite[Section 3.3]{BO10}): The construction is focused on
reducing the diameters of the boxes as fast as possible instead of keeping
them uniform on every level. Hence, one has to store an individual transfer
matrix for almost every box and an individual coupling matrix for every
admissible leaf of the block cluster tree.
Since approximating a block with less than $ k^{2} $ entries would therefore
be inefficient, the recursion of the algorithm for the cluster tree
construction stops when it reaches a box with not more than $ 2k $
characteristic points. On the finest mesh ($ n = 2097152 $), $ \widetilde{G}_{conv} $ requires 23328.1 MB for the leaf matrices, 33718.1 MB for the transfer matrices, 241218.8 MB for the nearfield matrices and 206819.6 MB for the coupling matrices. $ \widetilde{K}_{conv} $ requires 17496.1 MB for the leaf matrices, 25187.2 MB for the transfer matrices, 200100.7 MB for the nearfield matrices and 126874.0 MB for the coupling matrices. The resulting $ L^{2} $-error is $ 7.022_{-4} $.

In order to test our translation-invariant compression method on a more challenging boundary containing edges and corners, we approximate the Neumann values of the function $ u_{1} : \bbbr^{3} \rightarrow \bbbr, x \mapsto x_{1}^{2} - x_{3}^{2}$ on the boundary $ \Gamma_{C} := \partial([-1,1]\times[-\tfrac{3}{4}, \tfrac{3}{4}] \times [-\tfrac{1}{2} , \tfrac{1}{2}]) $ of a rectangular cuboid by the same method as before, i.e., we replace the matrices $ G $ and $ K $ appearing in (\ref{Galerkin_LS}) with their respective translation-invariant $ \htwo $-matrix approximations $ \widetilde{G} $ and $ \widetilde{K} $ and then solve the resulting perturbed linear system with the cg method. This time, we use $ r^{(0)} = [-1,1]\times[-\tfrac{3}{4}, \tfrac{3}{4}] \times [-\tfrac{1}{2} , \tfrac{1}{2}] $ as the starting box,
$ \const_{rk} = 1 $ for the determination of $ \ell_{max}^{I \times J} $ and increasing quadrature orders. The number of quadrature points per dimension for the regular integrals is denoted by $ q_{reg} $ and the number of quadrature points per dimension for the singular integrals is denoted by $ q_{sing} $. The characteristic points, $ \eta $, $ \ell_{max}^{I \times I} $ and the accuracy of the cg method are chosen as before.

\begin{table}[H]
	\centering
	\begin{tabular}{ |r|r|r|r|r|r|r|r| }
		$ n $ & $ q_{reg} $ & $ q_{sing} $  & $ \theta $ & $ k$ &  $ \ell_{max}^{I\times I} $ & $ \ell_{max}^{I\times J} $ &	$ \epsilon_{L^{2}} $ \\ \hline  
		   8112 & 4 &  6 & 4 & 125 &  9 &  9 & $ 3.895_{-2} $ \\
		  32448 & 4 &  7 & 5 & 216 & 12 & 12 & $ 1.484_{-2} $ \\
		 129792 & 4 &  8 & 6 & 343 & 12 & 12 & $ 7.408_{-3} $ \\
		 519168 & 5 &  9 & 7 & 512 & 15 & 15 & $ 2.626_{-3} $ \\
		2076672 & 5 & 10 & 8 & 729 & 18 & 18 & $ 1.332_{-3} $
	\end{tabular}
	\caption{Parameters and resulting $ L^{2} $-error $\epsilon_{L^{2}} := \|g_{N} - \tilde{g}_{N}\|_{L^{2}(\Gamma_{C})}$.}
	\label{par_err_tab_cuboid}
\end{table}

\begin{table}[H]
	\centering
	\begin{tabular}{ |r|r|r|r|r| }
		$ n $  & \small leaf matrices & \small transfer matrices & \small nearfield matrices & \small coupling matrices\\\hline  
		8112    &    15.5    &      4.3 &    55.9    &     48.0   \\
		32448   &    107.0   &     17.1 &    192.4   &     215.8  \\
		129792  &    679.4   &     43.1 &    1906.6  &     484.8  \\
		519168  &    4056.4  &    120.0 &    6915.7  &     1488.1 \\
		2076672 &    23101.6 &    291.9 &    26143.3 &     3843.8 
	\end{tabular}
	\caption{Storage requirements of $ \widetilde{G} $ on $ \Gamma_{C} $ in MB.}
	\label{storage_Gtilde_cuboid_tab}
\end{table}

\begin{figure}[H]
	\begin{tikzpicture}
	\begin{axis}[
	xmode = log,
	log basis x={2},
	ymode = log,
	log basis y={2},
	xtick distance = 4.0,  
	ytick distance = 4.0,
	grid = both,
	width = 0.5\textwidth,
	height = 0.5\textwidth,
	xlabel = {$ n $},
	ylabel = {\small Storage in MB},
	x label style={at={(axis description cs:0.5,0.02)}},
	y label style={at={(axis description cs:0.06,.5)}},
	legend pos = north west,
	legend style={nodes={scale=0.5}}]
	\newcommand\m{(log2(x)/2.0) - 2.4929209685};
	\newcommand\rank{(\m + 1)^3};
	\addplot[smooth, thick, red, domain = 8112:2076672]{16 * x *\rank * 8.0/(1024.0 * 1024.0)};
	\addplot[smooth, thick, teal, mark = x] table[x={n}, y={nearfield}] {Ktilde_cuboid.dtn};
	\addplot[smooth, thick, cyan, mark = cube] table[x={n}, y={leaf}] {Ktilde_cuboid.dtn};
	\legend{\large$ \mathcal{O}(nk) $, \large nearfield matrices, \large leaf matrices}
	\end{axis}
	
	\begin{axis}[xshift = 7.5cm,
	xmode = log,
	log basis x={2},
	ymode = log,
	log basis y={2},
	xtick distance = 4.0,  
	ytick distance = 4.0,
	grid = both,
	width = 0.5\textwidth,
	height = 0.5\textwidth,
	xlabel = {$ n $},
	ylabel = {\small Storage in MB},
	x label style={at={(axis description cs:0.5,0.02)}},
	y label style={at={(axis description cs:0.06,.5)}},
	legend pos = north west,
	legend style={nodes={scale=0.5}}]
	\newcommand\m{(log2(x)/2.0) - 2.4929209685};
	\newcommand\rank{(\m + 1)^3};
	\addplot[smooth, thick, red, domain = 8112:2076672]{64 * log2(x) *\rank^2 * 8.0/(1024.0 * 1024.0)};
	\addplot[smooth, thick, teal, mark = x] table[x={n}, y={coupling}] {Ktilde_cuboid.dtn};
	\addplot[smooth, thick, cyan, mark = cube] table[x={n}, y={transfer}] {Ktilde_cuboid.dtn};
	\legend{\large$ \mathcal{O}(\log_{2}(n)k^{2}) $, \large coupling matrices, \large transfer matrices}
	\end{axis}
	\end{tikzpicture}
	\caption{Storage requirements of $ \widetilde{K} $ on $ \Gamma_{C} $ in MB.}
	\label{storage_Ktilde_cuboid_fig}
\end{figure}
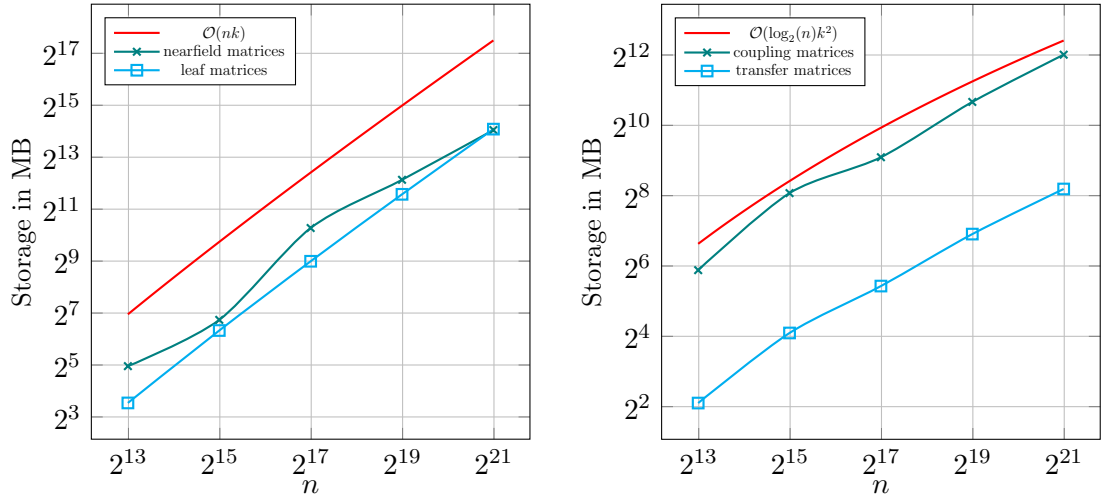

As in the previous experiment, the storage requirements of the translation-invariant $ \htwo $-matrix approximations and the resulting $ L^{2} $-error of the Neumann values behave as predicted by the theory (cf. Tables \ref{par_err_tab_cuboid}, \ref{storage_Gtilde_cuboid_tab}, \ref{storage_Ktilde_cuboid_tab} and Figures \ref{storage_Gtilde_cuboid_fig}, \ref{storage_Ktilde_cuboid_fig}).

\begin{table}[H]
	\centering
	\begin{tabular}{ |r|r|r|r|r| }
		$ n $  & \small leaf matrices & \small transfer matrices & \small nearfield matrices & \small coupling matrices \\ \hline
		 8112   & 11.6    &  4.3  & 31.0    & 59.0   \\
		32448   & 80.3    & 17.1  & 106.7   & 269.2  \\
		129792  & 509.6   & 43.1  & 1237.3  & 545.8  \\
		519168  & 3042.4  & 120.0 & 4480.0  & 1624.1 \\
		2076672 & 17326.6 & 291.9 & 16913.0 & 4119.6  
	\end{tabular}
	\caption{Storage requirements of $ \widetilde{K} $ on $ \Gamma_{C} $ in MB.}
	\label{storage_Ktilde_cuboid_tab}
\end{table}

\begin{figure}[H]
	\begin{tikzpicture}
	\begin{axis}[
	xmode = log,
	log basis x={2},
	ymode = log,
	log basis y={2},
	xtick distance = 4.0,  
	ytick distance = 4.0,
	grid = both,
	width = 0.5\textwidth,
	height = 0.5\textwidth,
	xlabel = {$ n $},
	ylabel = {\small Storage in MB},
	x label style={at={(axis description cs:0.5,0.02)}},
	y label style={at={(axis description cs:0.06,.5)}},
	legend pos = north west,
	legend style={nodes={scale=0.5}}]
	\newcommand\m{(log2(x)/2.0) - 2.4929209685};
	\newcommand\rank{(\m + 1)^3};
	\addplot[smooth, thick, red, domain = 8112:2076672]{16 * x *\rank * 8.0/(1024.0 * 1024.0)};
	\addplot[smooth, thick, teal, mark = x] table[x={n}, y={nearfield}] {Gtilde_cuboid.dtn};
	\addplot[smooth, thick, cyan, mark = cube] table[x={n}, y={leaf}] {Gtilde_cuboid.dtn};
	\legend{\large$ \mathcal{O}(nk) $, \large nearfield matrices, \large leaf matrices}
	\end{axis}
	
	\begin{axis}[xshift = 7.5cm,
	xmode = log,
	log basis x={2},
	ymode = log,
	log basis y={2},
	xtick distance = 4.0,  
	ytick distance = 4.0,
	grid = both,
	width = 0.5\textwidth,
	height = 0.5\textwidth,
	xlabel = {$ n $},
	ylabel = {\small Storage in MB},
	x label style={at={(axis description cs:0.5,0.02)}},
	y label style={at={(axis description cs:0.06,.5)}},
	legend pos = north west,
	legend style={nodes={scale=0.5}}]
	\newcommand\m{(log2(x)/2.0) - 2.4929209685};
	\newcommand\rank{(\m + 1)^3};
	\addplot[smooth, thick, red, domain = 8112:2076672]{64 * log2(x) *\rank^2 * 8.0/(1024.0 * 1024.0)};
	\addplot[smooth, thick, teal, mark = x] table[x={n}, y={coupling}] {Gtilde_cuboid.dtn};
	\addplot[smooth, thick, cyan, mark = cube] table[x={n}, y={transfer}] {Gtilde_cuboid.dtn};
	\legend{\large$ \mathcal{O}(\log_{2}(n)k^{2}) $, \large coupling matrices, \large transfer matrices}
	\end{axis}
	\end{tikzpicture}
	\caption{Storage requirements of $ \widetilde{G} $ on $ \Gamma_{C} $ in MB.}
	\label{storage_Gtilde_cuboid_fig}
\end{figure}
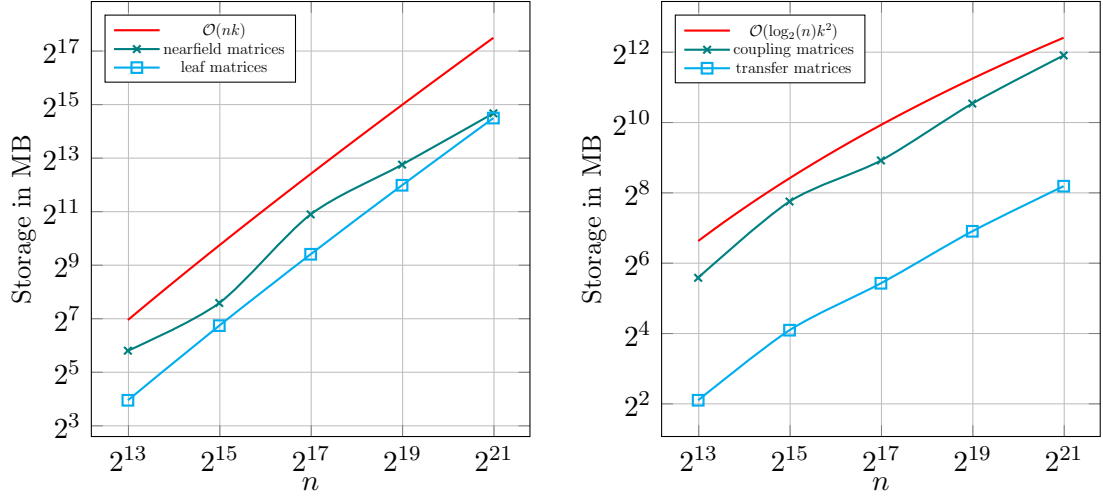

We again want to compare the results of the translation-invariant $ \htwo $-matrix compression with the corresponding results of the conventional $ \htwo $-matrix compression on the finest mesh ($ n = 2076672 $): 
$ \widetilde{G}_{conv} $ requires 23100.3 MB for the leaf matrices, 26630.6 MB for the transfer matrices, 247629.2 MB for the nearfield matrices and 167992.2 MB for the coupling matrices. $ \widetilde{K}_{conv} $ requires 17325.2 MB for the leaf matrices, 20321.7 MB for the transfer matrices, 193189.7 MB for the nearfield matrices and 115508.9 MB for the coupling matrices. The resulting $ L^{2} $-error is $ 1.503_{-3} $.

As the complexity bounds suggest and the results demonstrate, the
translation-invari\-ant approach reduces the storage requirements for the
transfer and coupling matrices drastically. Furthermore, the
translation-invariant $ \htwo $-matrix approximations require remarkably less
storage for the nearfield matrices than their respective conventional
counterparts. This is due to the fact that the leaves of the cluster trees
used in the translation-invariant approach contain significantly less
characteristic points than the leaves of the conventional cluster trees. The
resulting increase of the farfield can be handled very efficiently, since no
individual but only very few coupling matrices have to be stored.

We conclude that by properly exploiting the translation-invariance property
(\ref{eq:translation}) of the kernel function the storage
requirements of $ \htwo $-matrix approximations can be greatly reduced without affecting the desired accuracy.

\bibliographystyle{plain}
\bibliography{hmatrix}

\begin{thebibliography}{10}

\bibitem{AN92}
C.~R. Anderson.
\newblock An implementation of the fast multipole method without multipoles.
\newblock {\em SIAM J. Sci. Stat. Comp.}, 13:923--947, 1992.

\bibitem{BAHA05}
L.~Banjai and W.~Hackbusch.
\newblock {${\mathcal H}$}- and {${\mathcal H}^2$}-matrices for low and high
  frequency {Helmholtz} equations.
\newblock {\em IMA J. Numer. Anal.}, 28:46--79, 2008.

\bibitem{BE00a}
M.~Bebendorf.
\newblock Approximation of boundary element matrices.
\newblock {\em Numer. Math.}, 86(4):565--589, 2000.

\bibitem{BEVE12}
M.~Bebendorf and R.~Venn.
\newblock Constructing nested bases approximations from the entries of
  non-local operators.
\newblock {\em Numer. Math.}, 121(4):609--635, 2012.

\bibitem{BESC21}
T.~Betcke and M.~W. Scroggs.
\newblock Bempp-cl: A fast {Python} based just-in-time compiling boundary
  element library.
\newblock {\em J. Open Source Software}, 6(59), 2021.
\newblock available at \url{https://doi.org/10.21105/joss.02879}.

\bibitem{BO10}
S.~{B\"orm}.
\newblock {\em Efficient Numerical Methods for Non-local Operators: {${\mathcal
  H}^2$}-Matrix Compression, Algorithms and Analysis}, volume~14 of {\em EMS
  Tracts in Mathematics}.
\newblock EMS, 2010.

\bibitem{BOBO18}
S.~{B\"orm} and C.~{B\"orst}.
\newblock Hybrid matrix compression for high-frequency problems.
\newblock {\em SIAM J. Matrix Anal. Appl.}, 41(4):1704--1725, 2020.

\bibitem{BOCH14}
S.~{B\"orm} and S.~Christophersen.
\newblock Approximation of integral operators by {Green} quadrature and nested
  cross approximation.
\newblock {\em Numer. Math.}, 133(3):409--442, 2016.

\bibitem{BOGR04}
S.~{B\"orm} and L.~Grasedyck.
\newblock Hybrid cross approximation of integral operators.
\newblock {\em Numer. Math.}, 101:221--249, 2005.

\bibitem{BOHA02a}
S.~{B\"orm} and W.~Hackbusch.
\newblock {${\mathcal{H}}^2$}-matrix approximation of integral operators by
  interpolation.
\newblock {\em Appl. Numer. Math.}, 43:129--143, 2002.

\bibitem{BOHA03}
S.~{B\"orm} and W.~Hackbusch.
\newblock Approximation of boundary element operators by adaptive {${\mathcal
  H}^2$}-matrices.
\newblock {\em Found. Comp. Math.}, 312:58--75, 2004.

\bibitem{BOLOME02}
S.~{B\"orm}, M.~{L\"ohndorf}, and J.~M. Melenk.
\newblock Approximation of integral operators by variable-order interpolation.
\newblock {\em Numer. Math.}, 99(4):605--643, 2005.

\bibitem{BOME15}
S.~{B\"orm} and J.~M. Melenk.
\newblock Approximation of the high-frequency {Helmholtz} kernel by nested
  directional interpolation: error analysis.
\newblock {\em Numer. Math.}, 137(1):1--34, 2017.

\bibitem{CORAZO15}
E.~Corona, A.~Rahimian, and D.~Zorin.
\newblock A tensor-train accelerated solver for integral equations in complex
  geometries.
\newblock {\em J. Comp. Phys}, 334:145--169, 2015.

\bibitem{DAFAGRHASA01}
W.~Dahmen, B.~Faermann, I.~G. Graham, W.~Hackbusch, and S.~A. Sauter.
\newblock Inverse inequalities on non-quasiuniform meshes and applications to
  the mortar element method.
\newblock {\em Math. Comp.}, 73:1107--1138, 2004.

\bibitem{ERSA98}
S.~Erichsen and S.~A. Sauter.
\newblock Efficient automatic quadrature in 3-d {Galerkin} {BEM}.
\newblock {\em Comput. Meth. Appl. Mech. Eng.}, 157:215--224, 1998.

\bibitem{DAFO09}
W.~Fong and E.~Darve.
\newblock The black-box fast multipole method.
\newblock {\em J. Comp. Phys.}, 228:8712--8725, 2009.

\bibitem{GI01}
K.~Giebermann.
\newblock {M}ultilevel approximation of boundary integral operators.
\newblock {\em Computing}, 67:183--207, 2001.

\bibitem{GIRO02}
Z.~Gimbutas and V.~Rokhlin.
\newblock A generalized fast multipole method for nonoscillatory kernels.
\newblock {\em SIAM J. Sci. Comput.}, 24(3):796--817, 2002.

\bibitem{GR04}
L.~Grasedyck.
\newblock Adaptive recompression of {${\mathcal H}$}-matrices for {BEM}.
\newblock {\em Computing}, 74(3):205--223, 2004.

\bibitem{GRRO87}
L.~Greengard and V.~Rokhlin.
\newblock A fast algorithm for particle simulations.
\newblock {\em J. Comp. Phys.}, 73:325--348, 1987.

\bibitem{GRRO97}
L.~Greengard and V.~Rokhlin.
\newblock A new version of the fast multipole method for the {Laplace} equation
  in three dimensions.
\newblock In {\em Acta Numerica 1997}, pages 229--269. Cambridge University
  Press, 1997.

\bibitem{HA92}
W.~Hackbusch.
\newblock {\em {E}lliptic {D}ifferential {E}quations. {T}heory and {N}umerical
  {T}reatment}.
\newblock Springer-Verlag Berlin, 1992.

\bibitem{HA95}
W.~Hackbusch.
\newblock {\em Integral equations}.
\newblock Birkh{\"{a}}user Basel, 1995.

\bibitem{HA15}
W.~Hackbusch.
\newblock {\em Hierarchical Matrices: Algorithms and Analysis}.
\newblock Springer, 2015.

\bibitem{HANO89}
W.~Hackbusch and Z.~P. Nowak.
\newblock On the fast matrix multiplication in the boundary element method by
  panel clustering.
\newblock {\em Numer. Math.}, 54(4):463--491, 1989.

\bibitem{HSWE08}
G.~C. Hsiao and W.~L. Wendland.
\newblock {\em Boundary Integral Equations}.
\newblock Number 164 in Appl. Math. Sci. Springer, 2008.

\bibitem{RO85}
V.~Rokhlin.
\newblock Rapid solution of integral equations of classical potential theory.
\newblock {\em J. Comp. Phys.}, 60:187--207, 1985.

\bibitem{SA96}
S.~A. Sauter.
\newblock Cubature techniques for 3-d {G}alerkin {BEM}.
\newblock In W.~Hackbusch and G.~Wittum, editors, {\em Boundary Elements:
  Implementation and Analysis of Advanced Algorithms}, pages 29--44.
  Vieweg-Verlag, 1996.

\bibitem{SASC11}
S.~A. Sauter and C.~Schwab.
\newblock {\em Boundary Element Methods}.
\newblock Springer, 2011.

\bibitem{TY96}
E.~E. Tyrtyshnikov.
\newblock Mosaic-skeleton approximation.
\newblock {\em {C}alcolo}, 33:47--57, 1996.

\bibitem{TY99}
E.~E. Tyrtyshnikov.
\newblock Incomplete cross approximation in the mosaic-skeleton method.
\newblock {\em Computing}, 64:367--380, 2000.

\bibitem{BIYIZO04}
L.~Ying, G.~Biros, and D.~Zorin.
\newblock A kernel-independent adaptive fast multipole algorithm in two and
  three dimensions.
\newblock {\em J. Comp. Phys.}, 196(2):591--626, 2004.

\end{thebibliography}

\end{document}